\documentclass[twoside,a4paper,12pt,centertags]{amsart}
\usepackage{amsmath,amssymb,verbatim,vmargin}
\usepackage[bookmarks=true]{hyperref}   
\theoremstyle{plain}
\newtheorem{theorem}{Theorem}[section]
\newtheorem{lemma}[theorem]{Lemma}
\newtheorem{corollary}[theorem]{Corollary}
\newtheorem{proposition}[theorem]{Proposition}

\theoremstyle{definition}

\newtheorem{remark}[theorem]{Remark}

\numberwithin{equation}{section}
\title[{Flat forms, bi-Lipschitz parametrizations,
and smoothability}]{{Flat forms, bi-Lipschitz parametrizations,
and smoothability of manifolds}}
\author{Juha Heinonen
 and Stephen Keith}
\thanks{The first author was supported by the
NSF grants DMS 0244421, 0353549, and 0652915.}
\address{(J.H.) Department of Mathematics\\
University of Michigan\\
Ann Arbor, MI 48109, USA}
\thanks{
The second author was employed at the University of Helsinki and at
the Centre for Mathematics and its Application at Australia National University during
part of the research for this paper.}
\date{September 14, 2009}

\DeclareMathOperator{\dist}{dist}

\newcommand{\barint}{\mbox{$ave \int$}}

\def\barint_#1{\mathchoice
	  {\mathop{\vrule width 6pt height 3 pt depth -2.5pt
		  \kern -8pt \intop}\nolimits_{#1}}%
	  {\mathop{\vrule width 5pt height 3 pt depth -2.6pt
		  \kern -6pt \intop}\nolimits_{#1}}%
	  {\mathop{\vrule width 5pt height 3 pt depth -2.6pt
		  \kern -6pt \intop}\nolimits_{#1}}%
	  {\mathop{\vrule width 5pt height 3 pt depth -2.6pt
	\kern -6pt \intop}\nolimits_{#1}}}

\newlength{\defbaselineskip}

\newcommand{\Z}{\mathbf{Z}}             
\newcommand{\N}{\mathbf{N}}
\newcommand{\R}{\mathbf{R}}

\newcommand{\rn}{\R^{n}}

\newcommand{\cH}{\mathcal H}

\newcommand{\Om}{\Omega}
\newcommand{\vp}{\varphi}

\renewcommand{\bmod}{{\rm mod}}

\newcommand{\al}{\alpha}

\newcommand{\Va}{V_\alpha}

\DeclareMathOperator{\leng}{length}

\begin{document}
\bibliographystyle{acm}

\begin{abstract}
We give a sufficient condition for a metric (homology) manifold to be
locally bi-Lipschitz equivalent to an open subset in $\rn$. The condition
is a Sobolev condition for a measurable coframe of flat 1-forms.
In combination with an earlier work of D. Sullivan, our methods also yield an analytic characterization for smoothability of a Lipschitz manifold
in terms of a Sobolev regularity for frames in a cotangent structure.
In the proofs, we exploit  the duality
between flat chains and flat forms, and recently established
differential analysis on metric measure spaces. When specialized to $\rn$, our result gives a kind of asymptotic and Lipschitz version of the measurable Riemann mapping theorem as suggested
by Sullivan.

\end{abstract}

\maketitle

\section{Introduction.}\label{introduction}
On every topological manifold outside dimension four there is an
analytically defined complex of differential forms
that can be used, for example,
to develop Hodge theory, signature operators, and characteristic classes on the manifold.
This is classical in dimensions less than four, and can be done
smoothly \cite{moise:book}. The possibility of doing the same in dimensions above four (\cite{tel:ihes}, \cite{cst:qc}) comes from
two sources.
First,  we have Whitney's theory of Lipschitz invariant flat differential forms
that can be considered on every Lipschitz manifold \cite{wh:book}.
Second, we have
Sullivan's theorem which guarantees,
in said dimensions, a unique Lipschitz structure on every topological manifold \cite{sul:homeo}. There are well known obstructions
to smoothability of topological manifolds, so the Whitney theory cannot
always be done in a smooth framework \cite{munk:illi}, \cite{hima:smooth}, \cite{ks:essays}.

In the mid 1990s, invoking results from
geometric analysis, Sullivan proposed an idea how the Whitney flat
forms could be used to  detect whether a given Lipschitz manifold
possesses a smooth structure  \cite{sul:bld}, \cite{sul:found}, \cite{sul:video}.  He gave a definition for a cotangent structure over a Lipschitz manifold -- the definition involves a Lipschitz vector bundle over the manifold together with an identification
of its sections with flat forms on the manifold -- and then asked if the existence of such a structure implied that the manifold is smoothable.
An obstruction to smoothability emerged
in the form of a nontrivial local degree of a Lipschitz branched covering
map from the manifold to Euclidean space.

Flat forms can be considered on spaces more general than
Lipschitz manifolds. This follows from  the fundamental duality between
flat forms and flat chains discovered by Wolfe in 1948 \cite{wh:book}. Exploiting this fact together with recent differential analysis on singular
spaces, Sullivan and the first named author studied conditions that
would detect which metric spaces are locally bi-Lipschitz equivalent to Euclidean spaces \cite{hs:duke}. In addition to  natural topological
and geometric measure theoretic conditions on a given space, the
existence of a collection of appropriate flat forms
was stipulated. Such a collection, called a Cartan-Whitney presentation in \cite{hs:duke}, is analogous to a measurable coframe
on a Lipschitz manifold.
Here, too,
 an
obstruction to having  local bi-Lipschitz parametrizations emerged in the form of a nontrivial local degree of a Lipschitz branched covering map.

In this paper, we present  conditions
that remove the aforementioned obstructions; the conditions are similar  in both cases. In the first case, we have  the following result.

\begin{theorem}\label{sbility}
An oriented Lipschitz manifold
admits a smooth structure
if and only if it admits a cotangent structure in the sense of Sullivan  with local frames in the Sobolev space
$H^{1,2}$.
\end{theorem}

In Section \ref{sbilitysec}, we present a variant of Theorem \ref{sbility} more in line with the rest of this paper (Theorem \ref{sbility2}).

To the best of our knowledge, Theorems \ref{sbility} and \ref{sbility2} provide the first
analytic conditions for smoothability of manifolds; they can be viewed as  regularity results for coframes.  Recall that by Sullivan's work \cite{sul:homeo},
Theorem \ref{sbility} gives a necessary and sufficient condition
for a topological manifold outside dimension four to be smoothable. (As mentioned earlier, dimensions below four are of course classical.)
Note that Theorem \ref{sbility} has no dimensional restrictions.

Besides the results in  \cite{munk:illi}, \cite{hima:smooth}, \cite{sul:bld},
a geometric condition for smoothability of Lipschitz manifolds  was given by Whitehead in \cite{whi:trans}.
(See also \cite{pugh:smooth}.)

\

In the case of bi-Lipschitz parametrizations of metric spaces, we have the following result which
verifies a conjecture made in \cite[Remark 2.5]{hs:duke}. The required notions are explained
later in this introduction.

\begin{theorem}\label{main}
Let
$X\subset \R^N$ be a locally Ahlfors $n$-regular and
 locally linearly contractible homology $n$-manifold, $n\ge2$.
 If $X$ admits
local Cartan-Whitney presentations in the  Sobolev
space $H^{1,2}$, then $X$ admits
local  bi-Lipschitz parametrizations by
 $\R^n$. \end{theorem}

The conditions of $n$-manifold, local Ahlfors $n$-regularity, and local linear contractibility are obvious necessary conditions
for a space to admit local bi-Lipschitz parametrizations by
$\R^n$.  These conditions  are not sufficient, not even for subsets of $\R^N$ admitting
Cartan-Whitney presentations (except for $n=2$; see the next paragraph.)
 We refer to
 \cite{ss:complexes}, \cite{sem:nonexistence}, \cite{sem:parametrizations}, \cite{laa:blms}, \cite{bis:a1}, \cite{hs:duke}, \cite{hr:duke}, \cite{hei:icm} for examples and further discussion.

In dimension $n=2$,  the Sobolev condition for the Cartan-Whitney presentation turns out to be redundant. Local Ahlfors 2-regularity, local linear  contractibility, and the existence of local Cartan-Whitney
presentations characterize the 2-surfaces in Euclidean space
that admit local bi-Lipschitz parametrizations
by $\R^2$.
This observation (made by M. Bonk and the first named author;
see \cite[Theorem 2]{hei:icm})
uses the main result of \cite{hs:duke}, a local version of the quasisymmetric parametrization theorem
of Bonk and Kleiner \cite{bk:invent} proved by Wildrick  \cite{wild:thesis}, \cite{wild:stru}, and the  Morrey's measurable Riemann
mapping theorem.
All the aforementioned conditions including the existence of local Cartan-Whitney
presentations are required for this characterization, as follows from works by Semmes \cite{sem:bilipschitz} and  Bishop \cite{bis:a1}.

We note that it is a corollary to Theorem \ref{main} that a space satisfying the hypotheses
of the theorem is a topological manifold.
In this way, we arrive at an analytic characterization of topological manifolds among homology manifolds, quite different from the topological characterization due to Cannon, Edwards, Quinn, and others; see \cite{can:shr},  \cite{edw:icm}, \cite{quinn:mmj}, \cite{bfmw:gen}.
More precisely, we can state the following corollary.

\begin{corollary}\label{mancoro}
Let $X$ be a homology $n$-manifold, $n\ge2$. Then $X$ is a manifold if and only if every point in $X$ has a neighborhood that can be embedded in some $\R^N$ such that the image of the neighborhood is  locally Ahflors $n$-regular, locally linearly contractible, and admits local Cartan-Whitney presentations in the Sobolev space $H^{1,2}$.
\end{corollary}

 For each $n\ge3$, there are locally Ahlfors $n$-regular and
 locally linearly contractible subsets of some $\R^N$ that admit Cartan-Whitney presentations and are homology $n$-manifolds, but not manifolds; see \cite[9.1]{hr:duke}, \cite[2.4]{hs:duke}. Thus a (Sobolev) regularity condition as in Corollary \ref{mancoro} is necessary for $n\ge3$. We also point out that none of the exotic homology manifolds constructed in \cite{bfmw:gen} can serve as examples to this end, by the locality of the Quinn invariant \cite{quinn:mmj} and by \cite[Corollary 2.3]{hs:duke}.

There has been considerable interest in recent years
in trying to find intrinsic characterizations for classes of metric
spaces that admit local bi-Lipschitz parametrizations
by $\R^n$.
A simple geometric characterization of such spaces, especially in
dimensions $n\ge3$, seems to be a difficult problem.
However, a number of nontrivial sufficient conditions have emerged in the works
\cite{toro:jdg}, \cite{toro:duke},
\cite{ms:jdg}, \cite{bhr:crelle}, \cite{bhs:ab}, \cite{bl:bilip}, \cite{bhs:vuo}.
With the exception of \cite{toro:duke},   \cite{bhs:vuo}, these papers deal with the
case $n=2$. Theorem \ref{main} provides a new type of sufficient condition for a space to admit local bi-Lipschitz parametrizations by
$\R^n$ for all $n\ge2$.

Theorem \ref{main} follows from a more general result
(Theorem \ref{main2}), which we will formulate below.
First we require  some terminology and definitions.

\subsection{Basic definitions}\label{basic} A mapping $f:Y\to Z$
between metric spaces is  {\it bi-Lipschitz}
if there is a constant $L\ge1$ such that
$$
L^{-1}d_Y(a,b)\le d_Z(f(a),f(b))\le Ld_Y(a,b)
$$
for each pair of points $a,b\in X$.
We say that a metric
space  $X$ {\it admits
local  bi-Lipschitz parametrizations by}
 $\R^n$ if every point in $X$ has a neighborhood
that is bi-Lipschitz homeomorphic to an open
 subset of $\R^n$.
Throughout this paper, subsets of $\R^N$
are understood to carry the inherited metric.

A {\it topological $n$-manifold}, $n\ge2$,
is a separable and metrizable topological space that is locally homeomorphic
to open sets in $\rn$; if these homeomorphisms can be chosen so that the corresponding transition
functions are bi-Lipschitz, we have a {\it Lipschitz $n$-manifold}.
By \cite[4.2 and 4.5]{lv:elements}, every Lipschitz $n$-manifold
$X$ is (Lipschitz) homeomorphic to a subset  of $\R^{n(n+1)}$ that admits local bi-Lipschitz
parametrizations by $\rn$.
A Lipschitz manifold is {\it smoothable} if there is a further subsystem of charts that yield diffeomorphisms
as transition functions.

A (separable and metrizable) space $X$ is a {\it homology
$n$-manifold}, $n\ge2$, if $X$ is finite dimensional, locally compact and locally
contractible, and if the integral
homology groups satisfy
\begin{equation}\label{homoman}
H_*(X;X\setminus \{x\})=H_*(\R^n;\R^n\setminus\{0\})
\end{equation}
for each $x\in X$. If $X$ is an $n$-manifold, then
\eqref{homoman} holds, but for each $n\ge3$ there are nonmanifolds $X$
satisfying \eqref{homoman}.  See \cite{hs:duke} and \cite[Section 1]{hr:duke} for more discussion.

A metric space is said to be {\it locally Ahlfors $n$-regular}
if it has Hausdorff dimension $n$ and if for every compact
set $K$ in the space there exist numbers $r_K>0$ and $C_K\ge1$ such that
\begin{equation}\label{aregu}
C_K^{-1}r^n\le \cH_n(B(x,r))\le C_Kr^n
\end{equation}
for every metric ball $B(x,r)$ in the space, with center $x\in K$
and radius $r<r_K$. Here and later, $\cH_\alpha$, $\alpha>0$,
 denotes the $\alpha$-dimensional
Hausdorff measure.

A metric space is said to be {\it locally linearly contractible}
if for  every compact set
$K$ in the space there exist numbers $r_K>0$ and $C_K\ge1$ such that
every metric ball $B(x,r)$ in the space, with center
$x\in K$ and radius $r<r_K$, contracts to a point inside
$B(x,C_Kr)$.

For a locally Ahlfors $n$-regular and locally linearly contractible space $X$,
the dimension $n$ and the constants appearing in the preceding two conditions are called the {\it local data} of $X$ (typically near some point of interest in $X$).

Let $X\subset \R^N$ be a locally Ahlfors $n$-regular and
 locally linearly contractible homology $n$-manifold, $n\ge2$.
A local {\it Cartan-Whitney presentation}
on $X$  consists of an $n$-tuple $\rho=(\rho_1,\dots,\rho_n)$
of flat 1-forms in the sense Whitney \cite{wh:book}, defined
in an open $\R^N$-neighborhood of a point $p$ in $X$,
such that when restricted
to $X$ the associated ``volume form'' satisfies
\begin{equation}\label{vol}
{\mbox{essinf}}\, *(\rho_1\wedge\dots\wedge\rho_n)>0\,.
\end{equation}
The restriction makes sense under the present assumptions. See Sections 2 and 3 for a more precise discussion.

The existence of Cartan-Whitney presentations is a
bi-Lipschitz invariant condition, for flat forms can be pulled
back by Lipschitz maps.

The following should be considered the main theorem of this paper.

\begin{theorem}\label{main2}
Let $X\subset \R^N$ be a locally Ahlfors
 $n$-regular and
 locally linearly contractible homology $n$-manifold, $n\ge2$.
Assume that $X$ admits a
local Cartan-Whitney presentation $\rho=(\rho_1,\dots,\rho_n)$ in the  Sobolev
space $H^{1,2}$ near a point $p\in X$.
 Then the mapping
defined by
\begin{equation}\label{main2k}
F_p(x):=\int_{[p,x]}\rho\,,
\end{equation}
for $x\in \R^N$ sufficiently close to $p$, is bi-Lipschitz in a
neighborhood
of $p$ in $X$.
\end{theorem}

It is clear that Theorem \ref{main2} implies Theorem \ref{main}.

In \eqref{main2k}, the integration of the $n$-tuple of flat forms can be interpreted via the duality of flat chains and flat forms; this duality also has a more precise meaning as an integral as explained in \ref{intrinsic}.
The main fact   is that
by radially integrating flat forms satisfying \eqref{vol}, one obtains a Lipschitz branched
covering mapping, or a mapping of bounded
length distortion \cite{mv:bld}, \cite{hr:duke}.
 This important observation was made by Sullivan \cite{sul:bld}
for mappings in $\rn$ by using  fundamental results of Reshetnyak
\cite{resh:spacemaps},
\cite{resh:mappings}. A similar statement was proved for more singular spaces in \cite{hs:duke}.
The task in this paper is to show that the local degree of the
branched covering map $F_{p}$ as defined in \eqref{main2k} is one, under the given Sobolev regularity
assumption.

Next we consider the situation when $X=\R^n$. Theorem \ref{main2} is already of interest in this case.

\begin{theorem}\label{main3}
Let
$\Omega\subset \rn$, $n\ge2$, be an open set and let $\rho\in H^{1,2}(\Om)$ be a Cartan-Whitney presentation.
Then for every $p\in\Om$ there exists $r_0>0$
and a bi-Lipschitz mapping
$F_p:B(p,r_0)\to\rn$ such that
\begin{equation}\label{main3k}
    ||\rho-dF_p||_{\infty, B(p,r)}\,\le\,
r\,||d\rho||_{\infty, B(p,r)}
\end{equation}
for all $r<r_0$, where $B(p,r)$ denotes an open ball of radius $r$ centered at $p$ and
$||\cdot||_{\infty, B(p,r)}$ stands for
the $L^{\infty}$-norm in $B(p,r)$.
\end{theorem}

Theorem \ref{main3} follows from Theorem \ref{main2} and from the proof of Theorem 4.2 in \cite[p.\ 30]{hs:duke} --  cf.\  equation \eqref{asymp} below.

The assumption $\rho\in H^{1,2}$ in Theorem \ref{main3} can be viewed as an
{\it integrability condition} sought by Sullivan  in \cite{sul:video}. Fashioning Sullivan's  terminology, we can rephrase the assertion in
Theorem \ref{main3} by saying that
{\it{every Cartan-Whitney presentation in the Sobolev space $H^{1,2}$ is asymptotically locally standard.}} In this way,
Theorem \ref{main3} can be seen as a form
of the measurable Riemann mapping theorem in the Lipschitz context.
We refer to \cite{sul:video} for an interesting further discussion; see also \cite[Section 6]{hei:liplec}.

To demonstrate the sharpness of Theorems \ref{main2} and \ref{main3}, consider
the case where $X=\R^n$, $n\ge2$,
and $\rho=dF$ for the map
\begin{equation}\label{wmap}
F(r,\theta,z)=(r,2\theta, z)
\end{equation}
 in
cylindrical coordinates. Then $\rho$ determines
a Cartan-Whitney presentation in $\R^n$,
$\rho$ belongs to
 $H^{1,2-\epsilon}$ for each $\epsilon>0$ in a fixed neighborhood
of $0\in\rn$,  and $F=F_0$ as in \eqref{main2k},
but $F$ is not bi-Lipschitz.

We have $d\rho=0$ in the preceding example. Indeed,
Theorem \ref{main2} was proved by Kilpel\"ainen and the first named author
in the case where $X$ is an open subset of $\R^n$ and
$d\rho=0$ \cite{hki:bld}.
We  use also here the key idea of \cite{hki:bld}
which combines the fact that a Sobolev function has
Lebesgue points in a large set with the fact that
the image of the branch set of a discrete and open mapping
has large Hausdorff dimension. In the present case, however, essential new difficulties arise from two sources:
the Whitney forms are not assumed to be closed and the
underlying set $X$ is not assumed to be smooth.
We overcome these difficulties
first by a more substantial use of the Whitney
theory \cite{wh:book} and by invoking some fundamental results from geometric measure theory.
The most crucial help comes from
  some rather deep
recent results and techniques from analysis on metric spaces \cite{sem:curves}, \cite{hk:quasi}, \cite{ch:lipschitz},
\cite{keith:indi}.

A Cartan-Whitney presentation is a kind of measurable coframe.
Theorem \ref{main} asserts that if such coframes exist
locally on $X$ with certain degree of differential regularity, then $X$ admits local bi-Lipschitz
parametrizations by $\R^n$. The required membership in
a Sobolev space is second order information, akin to curvature.
The precise definition in the present context is
 given in Section \ref{sobolev}.
We will point out in Remark \ref{vmormks} that in all our main results
this second order information can be replaced by
a VMO or vanishing mean oscillation type condition, with same conclusions. In fact, the
proof of Theorem \ref{main} is a reduction of the Sobolev
condition to a VMO condition.

Finally, we remark that the questions addressed in this paper make sense also
in dimension $n=1$, but the answers are easy and well known in this case. Thus, we assume from now on that $n\ge2$ unless otherwise specifically mentioned.

The paper is organized as follows. Sections \ref{sec2} and \ref{sec3} contain preliminary material. Section \ref{comp} contains a discussion of tangent cones and compactness in the present context.
Sections \ref{lebsection} and \ref{tmaps} contain the most technical parts of the proof of Theorem \ref{main2}, which is completed in Section \ref{pofmain}.  In the last section, Section \ref{sbilitysec}, we prove the  smoothability  theorem \ref{sbility} and its variant Theorem \ref{sbility2}.

\

{\bf Acknowledgments.}  The authors wish to thank Mario Bonk, Seppo Rickman, Karen E.\ Smith, and Dennis Sullivan for useful discussions. The first named author in particular wishes to express his sincere gratitude to Dennis Sullivan who taught him so much.

\
\

\emph{Juha Heinonen passed away during the very final preparation of this paper. S.K. dedicates this work to his memory.}

\

\section{Geometric measure theory.}\label{sec2}
In this section, we discuss the measure theoretic implications
of the topological and geometric conditions appearing
in the introduction.
The standard references for the ensuing
concepts and facts from geometric measure theory are \cite{fed:gmt}, \cite{mat:geometry}.

\subsection{Rectifiable sets.}\label{basicdef}
A metric space  is said to be {\it $n$-rectifiable}
if
it can be covered, up to a set of Hausdorff
$n$-measure zero,  by a countable number of Lipschitz images
of subsets of  $\R^n$. In particular, $n$-rectifiable sets have Hausdorff
dimension $n$.

We call a metric space {\it metrically $n$-dimensional}
if it is $n$-rectifiable and locally Ahlfors
$n$-regular.
From now on, in
connection with metrically $n$-dimensional spaces, phrases
such as  ``almost everywhere''
always refer
to the Hausdorff $n$-measure.

An Ahlfors $n$-regular subset of $\R^N$  need not be $n$-rectifiable, but it turns out that under the (quantitative) topological hypotheses of Theorem
\ref{main}, this will be the case. The following proposition may well
be known to the experts, but we have not seen it in print. (The proposition is also valid for $n=1$, in which case it is
well
known.)

\begin{proposition}\label{sech}
Every locally linearly contractible and
locally Ahlfors $n$-regular metric homology $n$-manifold $X\subset \R^{N}$ is metrically $n$-dimensional.
     \end{proposition}

\begin{proof}
Under the current hypotheses, it follows from Semmes's work \cite{sem:curves}
that  $X$ admits (locally) a
Poincar\'e inequality (cf.\ Section \ref{sobolev} below). This
understood, a theorem of Cheeger \cite[Theorem 14.2]{ch:lipschitz}
 implies that $X$ is $n$-rectifiable.
\end{proof}

\subsection{Measurable tangent bundles.}\label{tang}
If $X$ is a metrically
$n$-dimensional subset of $\R^N$, then
it possesses tangent $n$-planes at almost every point; the local
Ahlfors regularity guarantees that the {\it approximate
tangent $n$-planes} are genuine tangent planes \cite[3.2]{fed:gmt}.
For our purposes, the tangent planes are best
described in terms of a {\it Hausdorff convergence}.
 Thus, a set $X\subset \R^N$ admits
a tangent $n$-plane at a point $x\in X$, if
the sets
\begin{equation}\label{xr}
X_{x,r}:=r^{-1}(X-x)=\{y\in \R^N: ry+x\in X\}\,, \ \ \ r>0\,,
\end{equation}
converge to an $n$-dimensional vector subspace
of $\R^N$ locally in the Hausdorff distance as $r\to 0$.
The limit vector subspace, whenever it exists,
 is called the {\it tangent $n$-plane
to $X$ at $x$}, and  denoted
by $T_xX$. (Compare with the discussion later in Section \ref{tangcones}.)

A sequence of sets $(F_i)$ in $\R^N$
is said to {\it converge
 locally in the Hausdorff distance} to a set $F\subset \R^N$ if
for each compact set $K\subset \R^N$  the sequence
$(F_i\cap K)$  converges in the Hausdorff distance
to $F\cap K$; see  \cite[2.10.21]{fed:gmt}, \cite[Chapter 8]{ds:fractals}.

If $U$ is an open subset of $X$, then
$T_xU=T_xX$. If $X$ is metrically $n$-dimensional,
then the collection of all tangent $n$-planes is
denoted by $TX$, and called (with some abuse of language)
the {\it measurable} {\it tangent
bundle} of $X$. We regard $TX$ as a subset of the Grassmann  manifold $G(N,n)$
of all $n$-dimensional vector subspaces of $\R^N$.

Given a vector space
$V$, we let  $\wedge_kV$
and $\wedge^kV$ denote, respectively,
the  space of $k$-vectors and $k$-covectors on $V$.
An {\it orientation } of an $n$-plane $V\in G(N,n)$ is a choice
of a component of the one-dimensional subspace $\wedge_nV$
of $\wedge_n\R^N$. A measurable choice of orientations
for the tangent $n$-planes of a metrically
$n$-dimensional set is called an {\it orientation}
of its measurable tangent bundle. In other words, an orientation of $TX$
is an equivalence class of measurable functions $$
\zeta:X\to \wedge_n\R^N\,
$$
such that $\zeta(x)\in\wedge_nT_xX$ for almost every $x\in X$, where two functions $\zeta_1$ and
$\zeta_2$ are
declared equivalent provided there is a positive measurable function $c:X\to (0,\infty)$
satisfying  $\zeta_1(x)=c(x)\,\zeta_2(x)$ for almost
 every $x\in X$. By using the standard metric of $\R^N$, we have a canonical representative for each such equivalence class,
\begin{equation}\label{orient}
\xi:X\to \wedge_n\R^N\,, \ \
\ \ \  \ \ \xi(x)=v_1(x)\wedge\dots\wedge v_n(x)\,, \end{equation}
where $\xi(x)$ is the unit $n$-vector associated with
a collection of  $n$ linearly independent vectors $v_i(x)\in T_xX\subset
\R^N$
of unit length. See  \cite[1.6.1 and 3.2.25]{fed:gmt}.

Unless otherwise specifically stated, $\R^n$ is assumed to have its standard orientation $e_1\wedge\dots\wedge e_n$ determined by the standard orthonormal basis.

\subsection{Metric orientation.}\label{metor}
Let $X\subset \R^N$ be a homology $n$-manifold, as well as
metrically $n$-dimensional. Besides an orientation $\xi$ as
in \eqref{orient}, there is another orientation defined locally
from the topological data. Namely, every point in $X$ has
a neighborhood $U$ such that the compactly supported
(integral Alexander-Spanier) cohomology group satisfies
$H_c^n(U)=\Z$, and an {\it orientation} of $U$ is a choice
of a generator $g_U$ for this group.
Every connected subset $V$ of $U$ canonically inherits
an orientation from $U$, because
the inclusion $V\to U$ induces an isomorphism in cohomology.

It is important for our purposes to have an orientation
$\xi$ on $TX$ that locally matches the given topological orientation.
If $X$ is in addition
locally linearly contractible, then
such an orientation can be found \cite[3.10]{hs:duke}. More precisely,
let $U\subset X$ be open such that $H_c^n(U)=\Z$, and
choose a generator $g_U$ for $H_c^n(U)$. Then,
for almost every $x\in U$, the projection
$\pi_x:\R^N\to x+T_xX$ induces an isomorphism
$H_c^n(T_x)\to H_c^n(U)$, and
$\xi(x)$ is chosen so as to correspond to $g_U$ under this isomorphism.
(Note that the orientation $\xi(x)
\in\wedge_n\R^N$ canonically corresponds to an orientation in $H_c^n(T_x)$;
see, for example, \cite[Chapter 9]{milnor:char}.)
The procedure is explained in detail in \cite[3.4]{hs:duke};
see also \cite[5.9]{hr:duke}.

A pair of orientations $(g_U,\xi)$ on  $U$
for which the above compatibility condition holds
is said to be a {\it metric orientation} of $U$,
and $U$ is said to be {\it metrically oriented}
if such a pair is fixed.

\subsection{Duality and norms}\label{norms}
The standard Euclidean metric in $\R^N$ canonically determines an inner product structure, and hence a norm,
 in each $\wedge^k\R^N$ and $\wedge_k\R^N$. In particular, there is a canonical {\it Euclidean duality}
between the two spaces.
We will use an equivalent {\it comass} norm, and its dual {\it mass} norm, for covectors and vectors, respectively. Following \cite[1.8]{fed:gmt},  we denote these latter norms  by $||\cdot||$. Both comass and mass agree with the Euclidean norm for simple covectors and vectors.
Given an $n$-tuple of 1-covectors $\alpha=(\alpha_1,\dots,\alpha_n)\in(\wedge^1\R^N)^n$,
we write
\begin{equation}\label{tuplenorm}
||\alpha||:=(\sum_{i=1}^n||\alpha_i||^2)^{1/2}
\end{equation}
and
\begin{equation}\label{tuplenorm2}
\wedge\alpha:=\alpha_1\wedge\dots\wedge\alpha_n\in\wedge^n\R^N\,.
\end{equation}
Similar conventions hold for tuples of vectors in $(\wedge_1\R^N)^n$.
From now on it is understood that the vector spaces $\wedge^k\R^N$, $\wedge_k\R^N$, and
their finite products, are equipped with the norms $||\cdot||$ as in the preceding.

\subsection{Sobolev spaces on rectifiable sets.} \label{sobolev}
Let  $X\subset \R^N$ be an $n$-rectifiable set of locally finite Hausdorff $n$-measure.
We can define the
Sobolev space $H^{1,2}(X)$ as the (abstract) norm completion of the space
of all Lipschitz
functions $\varphi:X\to \R$ equipped with the norm
\begin{equation}\label{sobo}
||\varphi||_{1,2}:=||\varphi||_{L^2(X)}+||D\varphi||_{L^2(X)}\,,
\end{equation}
where $D\varphi$ stands for the {\it approximate differential}
of $\varphi$ which exists at almost every $x\in X$ as a linear map from the approximate
tangent plane $T_xX$ to $\R$ \cite[3.2.19]{fed:gmt}.
The metric in $\R^N$ gives the operator norm for the approximate differential, and in the $L^2$-norm
we use the Hausdorff $n$-measure.

If we assume moreover that $X$ is a metrically $n$-dimensional and locally linearly contractible
homology $n$-manifold,
the
Sobolev spaces just defined are known to enjoy properties similar to those
in Euclidean spaces. In particular, in this case, $H^{1,2}(X)$
is a vector subspace of $L^2(X)$ and every
$u\in H^{1,2}(X)$ is approximately differentiable almost
everywhere; we also have (locally) the following {\it Poincar\'e
inequality}:
\begin{equation}\label{poin}
\barint_{B\cap X}|u-u_{B\cap X}|^2\,d\cH_n\,\le\, C\,({{\mbox
{diam}}}B)^2\barint_{(\lambda B)\cap X}|Du|^2\,d\cH_n\,,
\end{equation}
where the barred integral sign denotes the mean-value, \begin{equation}\label{mvdef}
u_E:=\barint_Eu \, d\cH_n\,
\end{equation}
for a Borel set $E\subset X$,
 $B=B(x,r)$ is a
sufficiently small ball in $\R^N$ centered at $X$, and $\lambda B=B(x,\lambda r)$.
  The constants $C, \lambda \ge1$
depend only on the local data of $X$.

We will also speak of a membership in $H^{1,2}(X)$ of functions with values in a finite dimensional
Banach space -- more specifically $(\wedge^1\R^N)^n$ in our case, cf.\ the end of subsection \ref{intrinsic}. Such a  membership is understood
componentwise and is independent of the particular norm used; equivalently,  we can consider directly  Banach-valued Sobolev functions as in \cite{hkst:banach}.
Finally, the local Sobolev spaces $H^{1,2}_{loc}(X)$ are defined in an obvious way.

We refer to \cite{sem:curves}, \cite{hk:quasi}, \cite{ch:lipschitz}, \cite{hkst:banach}, \cite{hs:duke}, \cite{hr:duke},  \cite{kos:pisa},
\cite{keith:advances}, \cite{keith:indi}, \cite{KeithZhong} for a more complete
discussion on matters related to the Poincar\'e inequality. See in particular \cite[Section 3]{keith:indi} for the
equivalence of the above definition for the Sobolev space with
those appearing elsewhere.

\section{Whitney flat forms.}\label{sec3}
In this section, we discuss  flat forms
and  Cartan-Whitney presentations as defined
in \cite{hs:duke}.
Details for  the ensuing facts
can be found in  \cite{wh:book}, \cite{fed:gmt}.

\subsection{Flat forms.}
A (\emph{Whitney}) {\it flat k-form} in an open set $O\subset
\R^N$ is a $k$-form $\omega$ with measurable coefficients
satisfying
$$
\omega, \, d\omega\,\in L^{\infty}(O)\,,
$$
where $d\omega$ is understood in the sense of distributions; thus, $d\omega$ is the unique
$(k+1)$-form with bounded measurable coefficients satisfying
$$
\int_O\omega\wedge d\varphi=(-1)^{k+1}\int_{O}d\omega\wedge \varphi
$$
for all smooth  $(n-k-1)$-forms $\varphi$ with compact support in $O$.
 In this paper, we exploit the
description of flat $k$-forms  as the dual space of the normed space of flat $k$-chains.  The latter space
is defined as the completion of polyhedral $k$-chains
$P$ in $O$ with respect to the {\it flat norm}
$$
|P|_{{{\flat}}}:=\inf\{|P-\partial R|+|R|\}\,,
$$
where the infimum is taken over all polyhedral $(k+1)$-chains
$R$ in $O$, and where $|\cdot|$ denotes the {\it mass}
of a polyhedral chain. It is further required that every flat chain
has  compact support in $O$.
(Our considerations are all local, so that $O$ could be assumed to be convex throughout,
cf. \cite[p.\ 231]{wh:book}.)
By a theorem
of Wolfe  \cite[Chapter IX, Theorem 7C]{wh:book}
(see also \cite[4.1.19]{fed:gmt}),
the dual space of this normed space, the space of {\it flat $k$-cochains},
can be identified with the space of flat $k$-forms
equipped with the norm
\begin{equation}\label{fnorm}
||\omega||_{{{\flat}}}:=\max\{||\omega||_\infty,||d\omega||_\infty\}\,.
\end{equation}

Here on the right hand side, we use the pointwise comass of a form in determining the $L^\infty$-norm, cf.\  \ref{norms}.

It follows that flat $k$-forms
have a well defined
``trace'',  or restriction, on oriented $k$-simplexes
in $\R^N$ through the dual action. We often denote this action suggestively by an integration. For example, we write
$$
\langle \omega, [x,y]\rangle=\int_{[x,y]}\omega\,
$$
for an oriented line segment $[x,y]$ in $\R^N$,
if $\omega$ is a flat 1-form  defined in a neighborhood of $[x,y]$.
With this notation, we also have the Stokes theorem
\begin{equation}\label{w2}
\int_Pd\omega=\int_{\partial P}\omega\,
 \end{equation}
for  oriented polyhedral chains, or more generally for all flat chains.

A $k$-dimensional rectifiable set
of finite $k$-measure together with an orientation on its (approximate) tangent planes determines
a flat $k$-chain. Thus, every flat $k$-form has a well defined restriction, or action, on such a set
 \cite[4.1.28]{fed:gmt}.
In our context, we will need this fact
for flat $n$-forms acting
on metrically $n$-dimensional, metrically oriented sets,
and for flat 1-forms acting on rectifiable 1-dimensional sets.

We also require the fact that flat forms are closed under the wedge product, and can be pulled back by Lipschitz maps. Thus, if $\omega$ and
$\eta$ are flat forms in $\R^N$, then so is $\omega\wedge \eta$. Moreover, if $h:\R^M\to \R^N$ is a Lipschitz map, then the pullback $h^*(\omega\wedge\eta)$ can be defined; it is a flat form of the same degree as $\omega\wedge\eta$
in $\R^M$. We also have that $h^*(\omega\wedge\eta)=h^*(\omega)\wedge h^*(\eta)$.
See \cite[X.11]{wh:book}.

\subsection{Cartan-Whitney presentations.}\label{cwsub}
Let $X\subset \R^N$ be a metrically
$n$-dimensional and locally linearly contractible
homology $n$-manifold.
We assume that a metric orientation $(g_U,\xi)$
is chosen for an open
neighborhood $U$ of a point $p\in X$.
Then a ({\it local}) {\it Cartan-Whitney presentation near}
$p$ consists of an $n$-tuple
$$
\rho=(\rho_1,\dots,\rho_n)
$$
of flat 1-forms defined in an $\R^N$-neighborhood $O$ of $p$,
$O\cap X\subset U$,
such that
\begin{equation}\label{cw}
*(\rho_1\wedge\dots\wedge\rho_n)\ge c>0\,
\end{equation}
almost everywhere on $O\cap X$, where $c$ is a constant.
The wedge product $\rho_1\wedge\dots\wedge\rho_n$ is a
flat $n$-form and hence has a well defined action on every
(relatively) open patch $V$ of $O\cap X$, and \eqref{cw} is
interpreted as
\begin{equation}\label{cw2}
\int_V\rho_1\wedge\dots\wedge\rho_n\ge c\,\cH_n(V)>0\,
\end{equation}
for every such $V$.
An alternative, more intrinsic description of \eqref{cw} will be
provided in subsection \ref{intrinsic}. It is indicated
there how flat forms have well defined restrictions on rectifiable
sets. The Hodge $*$ operator used in \eqref{cw} then makes sense
in the usual way.

\subsection{Cartan-Whitney presentations and BLD-mappings}
\label{cartBLD}

Assume that a local Cartan-Whitney presentation $\rho$ is given in
an open convex $\R^N$-neighborhood $O$ of a point in $X$, where
$X\subset \R^N$ is  a metrically $n$-dimensional, metrically oriented, and locally linearly
contractible homology $n$-manifold. For each $q \in O\cap X$, define a
map
\begin{equation}\label{omap}
F_q(x):=\int_{[q,x]}\rho\,, \ \ \ \ x\in O\,.
\end{equation}
By \cite[Theorem 4.2]{hs:duke}, every such point $q$
has a
neighborhood  $U\subset X$ such that
the map $F_q$ is a uniformly BLD-mapping of $U$ into $\R^n$.
A
mapping $f:U\to \R^n$ is said to be an $L$-BLD-{\it mapping},
$L\ge1$,
 if $f$ is a sensepreserving, open and discrete $L$-Lipschitz map
such that
\begin{equation}\label{blddef}
L^{-1}\leng(\gamma)\le \leng(f(\gamma))\le L\leng(\gamma)
\end{equation}
for every curve $\gamma$ in $U$. Recall that a map is called {\it
discrete} if the preimage of every point under the map is a
discrete set. By a {\it uniformly} BLD-mapping in the present context we mean a BLD-mapping
such  that the constant $L$ in \eqref{blddef} depends only on the local data of $X$ near $q$ and on the flat norm of
$\rho$.

The proof of Theorem 4.2 in \cite[p.\ 30]{hs:duke} moreover gives the following asymptotic estimate
\begin{equation}\label{asymp}
||\rho-dF_q||_{\infty, B(q,r)}\,\le\,r\,||d\rho||_\infty\,, \ \ \ \ \ \ 0<r<r(q)\,.
\end{equation}

The {\it branch set} $B_f$ of a discrete and open map $f:U\to\R^n$
is the closed set of points in $U$ where $f$ does not determine a
locally invertible map. Under the present assumptions, the
branch set of a BLD-mapping is of measure zero and, if nonempty, of
positive Hausdorff $(n-2)$-measure. For the first assertion, see
\cite[Theorem 6.6]{hr:duke}. The second assertion follows from
\cite[Proposition 6.1]{hr:duke} and from the fact that the
image of a nonempty branch set of every discrete and open map from
a homology $n$-manifold to $\R^n$ must have positive Hausdorff
$(n-2)$-measure \cite[III 5.3, p.\ 74]{rick:quasiregular}.

In the present context, every BLD-mapping is locally bi-Lipschitz
in the complement of its branch set \cite[Theorem 4.2]{hs:duke}.
This fact is crucial to us.

We refer to  \cite{mv:bld} and \cite{hr:duke} for the general
theory of BLD-mappings.

\subsection{Restrictions of flat forms on $X$.}\label{intrinsic}
Let $X\subset\R^N$ be an $m$-rectifiable set of locally finite Hausdorff $m$-measure, $m\ge1$, and
let $\omega$ be a flat 1-form defined
in some open set $O$ of $\R^N$ with $U:=O\cap X$ nonempty.
We will canonically associate with each such  $\omega$ a measurable
almost everywhere defined
map
\begin{equation}\label{rest}
\omega\lfloor U:U\to \wedge^1\R^N\,, \ \ \ x\mapsto \omega\lfloor U(x)\in
T_x^*U\,,
\end{equation}
where $T_x^*U:=\wedge^1T_xU$
is viewed as the subspace
of the vector space $\wedge^1\R^N$ characterized by the following
property: $\eta\in T_x^*U$ if and only if $\langle \eta, v\rangle=0$
for every $v\in \wedge_1\R^N=\R^N$ perpendicular to the approximate tangent $m$-plane $T_xU$.  In this way, we can view
$\omega\lfloor U$
as a {\it measurable section} of a subbundle $T^*U$ (the {\it measurable cotangent bundle} of $U$)
 of the cotangent bundle $T^*\R^N$.

To this end, recall first that almost all
of $U$ can be covered by pairwise disjoint compact sets that are bi-Lipschitz
images of pieces of $\R^m$ \cite[3.2.18]{fed:gmt}. Consider
a compact set $K\subset\R^m$ and a bi-Lipschitz map
$g:K\to g(K)\subset U$. By standard Lipschitz extension theorems, $g$ can be extended to a Lipschitz
map $G:\R^n\to \R^N$. The issue is local, so we may further assume that $O$ is a ball
and that $G(\R^n)\subset O$.
It follows from \cite[Theorem  9A, p.\ 303]{wh:book} that for $\cH_m$ almost every $x\in G(K)$ the flat form
$\omega$ determines an element of $\wedge^1T_xU$. Moreover, this element is independent of the choice of $g$, and we denote it by $\omega\lfloor U(x)$. It follows that $\omega\lfloor U$ can be defined almost everywhere in $U$ as in  \eqref{rest}.

Assume now that $X$ is a metrically $n$-dimensional and locally linearly contractible homology $n$-manifold with a metric orientation $(g_U, \xi)$  chosen for $U$ as in \ref{cwsub}. The preceding construction shows that each
Cartan-Whitney presentation $\rho$  in $O$ determines an $n$-tuple of maps $U\to \wedge^1\R^N$. It follows that in
equation \eqref{cw2} one could equivalently integrate the almost
everywhere defined wedge product
${\rho_1}\lfloor U\wedge\dots\wedge{\rho_n}\lfloor U$ over $U$. That is,
\begin{equation}\label{cw3}
\int_U\rho_1\wedge\dots\wedge\rho_n
=\int_U\langle{\rho_1}\lfloor U(x)\wedge\dots\wedge{\rho_n}\lfloor U(x)\,,
\xi(x)\rangle\,d\cH_n(x)\,.
\end{equation}

Since the restriction $\omega\lfloor U$
of a Whitney flat 1-form $\omega$ as in \eqref{rest}
can be thought of as a vector valued function on $U$ (with values in $\wedge^1\R^N$), we can speak of a membership of $\omega$
in a Sobolev space on $U$ as discussed in \ref{sobolev}.
In particular, we can speak of a membership  in a Sobolev space
of a Cartan-Whitney presentation $\rho=(\rho_1,
\dots,\rho_n)$. The hypotheses in our main theorems should be understood in this vein.

\section{Tangent cones and flat compactness}\label{comp}
In this section, we record two facts about tangent cones that are crucial later in the paper.

\subsection{Tangent cones}\label{tangcones}
A {\it tangent cone} of a locally compact set $Y\subset \R^N$
at $y\in Y$ is a closed set $Z\subset \R^N$ for which there exists a sequence $r_i\to 0$ such that $Y_{y,r_i}\to Z$
locally in the Hausdorff distance. (See \ref{tang} and \eqref{xr}.)
Tangent cones exist for every locally compact $Y$ at every point, although they need not be unique. If $Z$ is an $n$-dimensional vector subspace of $\R^N$, or an {\it $n$-plane} for short, we call it a {\it tangent $n$-plane} at $y$ (even if it is not unique).
We sometimes write $Z=Z_y$ to emphasize the point where the rescaling takes place, or $Z=X_\infty$
if the point is clear from the context.
We refer to  \cite[Chapter 8]{ds:fractals} for  a careful discussion about tangent cones.

We require the following information about tangent cones of metrically $n$-dimen\-sion\-al
and locally linearly contractible homology $n$-manifolds.

\begin{proposition}\label{bb}
Let $X$ be a  metrically $n$-dimensional and locally linearly contractible homology $n$-manifold,
and let $X_p$ be a tangent cone of $X$ at a point $p\in X$.
 Then $X_p$ is a homology $n$-manifold satisfying the following global versions
 of Ahlfors regularity and linear local contractibility: there exists a constant $C\ge1$, depending only on the local data of $X$  near $p\in X$, such that $B(x,r)\cap X_p$ contracts to a point in $B(x,Cr)\cap X_p$
 and that
 \begin{equation}\label{aregu2}
 C^{-1}r^n\,\le\,\mathcal H_n(B(x,r))\,\le\,Cr^n
 \end{equation}
 whenever $x\in X_p$ and $0<r<\infty$.
 In particular, $X_p$ is metrically $n$-dimensional.
\end{proposition}

\begin{proof}
The fact that tangent cones of Ahlfors regular spaces are Ahlfors regular (of the same dimension)
is discussed in \cite[p.\ 61 ff.]{ds:fractals}; the discussion there extends to the local version as well.
The fact that $X_p$ satisfies the global contractibility condition and is a homology $n$-manifold, goes back to Borsuk \cite{borsuk:metr} and Begle \cite{begle:conv}, cf.\ \cite{grove:inv}. In the literature, this result has been stated for compact  spaces only, but the present local version can be proved similarly.
Finally, the last statement follows from Proposition \ref{sech}.
\end{proof}

We will also require the following uniform retraction property for linearly locally connected spaces.

\begin{proposition}\label{retr}
Let $X\subset \R^N$, $p\in X$, and $X_p$  be as in Proposition \ref{bb}. Assume that
$X_{p,r_i}\to X_p$. Then for every $R>0$ there exists $i_0$ with the following property: for every $i_0\le i,j\le \infty$
there are maps
\begin{equation}\label{retrk1}
\psi_{i,j}:B(p,R)\cap X_{p,r_i}\to B(p,2R)\cap X_{p, r_{j}}
\end{equation}
satisfying
\begin{equation}\label{retrk2}
|\psi_{i,j}(x)-x|\,\le\,C\,\dist(x, X_{p, r_{j}})\,,
\end{equation}
where $X_{p, r_\infty}=X_p$ and $C>0$ depends only on the local data of $X$ near $p$.
Moreover, the maps $\psi_{i,j}$ can be defined in all of $B(p,R)$ such that \eqref{retrk2} holds.
\end{proposition}

Proposition \ref{retr} follows from the local linear contractibility property in a customary way; see \cite[Section 5]{sem:curves}. In fact, the existence of retractions as in \eqref{retrk1}
are needed in the proof of Proposition \ref{bb}.

\subsection{Flat compactness} Let
 $X\subset \R^N$ be a  metrically $n$-dimensional
and locally linearly contractible homology $n$-manifold.
Assume that $U\subset X$ is open, precompact, and metrically oriented by $(g_U,\xi)$, as explained in \ref{metor}.
The pair $(U,\xi)$ defines a rectifiable $n$-current, hence a
flat $n$-chain \cite[4.1.28]{fed:gmt}.
If $x\in U$, the sets
$U_{x,r}$, as defined in \eqref{xr},
inherit the orientation from $U$ in a natural way, and,
therefore, we have a sequence $(U_{x,r})$ of flat $n$-chains.
We are
 interested in the compactness   of this
sequence  when $r$ is small.

It was proved in \cite[3.28]{hs:duke}
that $X$ is locally a cycle in the sense that every point
in $X$ has a neighborhood  such that $X$ restricted
to that neigborhood has no boundary in the sense of currents; that is,
\begin{equation}\label{noreuna}
\langle{\partial X,\omega\rangle}=\langle{X,d\omega\rangle}=0
\end{equation}
for each smooth $(n-1)$-form $\omega$ with support in the neighborhood.
We may clearly assume that $U$ is small enough so that \eqref{noreuna}
holds for all $\omega$ whose support intersected with $X$ lies in $U$.

Under these hypotheses,
we have the following result.

\begin{proposition} \label{prop:locr}
Given $x\in U$ and a sequence $(r_i)$ of positive real numbers tending to zero, we
can pass to a subsequence $(r_{i_j})$ so as to obtain a locally rectifiable
closed $n$-current as the limit of $(U_{x,r_{i_j}})$ in the sense of currents, as $r_{i_j}\to0$.
\end{proposition}

\begin{proof}
The compactness statement of this proposition follows essentially from the
Compactness Theorem \cite[4.2.17]{fed:gmt}. To conform with the
notation and assumptions of \cite{fed:gmt}, we need to use
\eqref{noreuna} and dispense with the boundary altogether. Let $B(0,s)$ denote an open ball in $\R^N$, centered at the
origin with radius $s>0$. We then have a  uniform mass bound guaranteed by
the assumptions,
\begin{equation}\label{massbound}
\cH_n(B(0,s)\cap U_{x,r})\le C r^n\,, \end{equation}
for all $s>0$, and for all small enough $r>0$,
where $C>0$ is independent of $r$. Write $U_i:=U_{x,r_i}$, and fix $s>0$.
We change each $n$-current $U_i$
 outside the ball $B(0,\frac{3}{2}s)$ as follows. Fix a smooth map
$\psi:\R^N\to\R^{N+1}$ such that $\psi|B(0,2s)$ is an
embedding, that $\psi|B(0,\frac{3}{2}s)=$identity, and that
$\psi(\R^N\setminus B(0,2s))=e_{N+1}=(0,\dots,0,1)\in\R^{N+1}$. Then for all sufficiently
 large $i$, the pushforward current $\psi_{\#}U_i$ satisfies
(cf. \cite[4.1.20]{fed:gmt})
\begin{equation}\label{push}
\partial \psi_{\#}U_i=0\,,  \ \ \ \ \psi_{\#}U_i|B(0,\frac{3}{2}s)=U_i\,.
\end{equation}
By using  \eqref{massbound}, \eqref{push}, and the
 Compactness Theorem \cite[4.2.17]{fed:gmt}, we obtain that the sequence $(\psi_{\#} U_i)$ is compact in the flat norm.
This understood, the proposition  follows by applying a
diagonalization argument in $i$ and $s$, by observing that local
convergence in the flat norm implies convergence in the sense of
currents, and by recalling that, as a consequence of the Closure Theorem
\cite[4.2.16]{fed:gmt}, (integral) flat chains with locally finite mass
are locally rectifiable. The proof is complete.
\end{proof}

\section{Tangents cones at  points of vanishing mean oscillation}\label{lebsection}
In this section,
we study the tangent cones
at a point of vanishing
mean oscillation for a Cartan-Whitney presentation.

\subsection{Basic set-up} \label{setup}Assume that a local
Cartan-Whitney presentation $\rho$ is given in a convex
$\R^N$-neighborhood $O$ of a point $p \in X$, where $X\subset \R^N$
is a metrically $n$-dimensional and locally linearly contractible
homology $n$-manifold (cf.\ Proposition \ref{sech}).  {\it It is not assumed at this point that $\rho$ belongs to
some Sobolev class.} Then fix a metrically oriented
neighborhood $U\subset X\cap O$ of $p$, as explained in \ref{metor}. For simplicity of notation,
we assume that $U=X\subset O$. Open balls in $\R^N$ are denoted by $B(x,r)$.

The preceding set-up understood, we simplify notation and write $\rho=(\rho_1,\dots,\rho_n)$ for the
restrictions $\rho\lfloor U=({\rho_1}\lfloor U, \dots, {\rho_n}\lfloor U)$ as defined in \ref{intrinsic}.
If $Z\subset X$ is a Borel set, we write
$$
\rho_Z=\barint_Z \rho \, d\cH_n
$$
for the element in $(\wedge^1\R^N)^n$, whose components are the mean values
 of the components of $\beta=(\beta_1,\dots,\beta_n)$
as usual Banach valued integrals, cf.\ \eqref{mvdef}.
Finally, recall the notation from \eqref{tuplenorm}, \eqref{tuplenorm2}.

Certain claims in the upcoming arguments are
said to hold everywhere in $X$, although strictly speaking one should consider
a fixed neighborhood of $p$, where the local assumptions on Ahlfors
regularity and linear contractibility hold, with fixed constants, together
with the concomitant analytic consequences such as the Poincar\'e inequality.
In  rescaling and asymptotic
arguments, such a localization becomes immaterial, and a small abuse of language
keeps the presentation simpler. Recall that the local data refers to the constants
appearing in the said local conditions; these include the flat norm of $\rho$ and the lower bound $c$ in \eqref{cw}.

\subsection{VMO-points and tangent cones}
A point $q\in X$ is said to be a {\it point of vanishing mean oscillation}, or a {\it VMO-point}, for $\rho$
if
\begin{equation}\label{vmo}
\lim_{r\to0}\barint_{B(q,r) \cap X}
||\rho-\rho_{B(q,r)\cap X}||
 \, d\cH_n=0\,.
\end{equation}

The following is the main result of this section.

\begin{proposition}\label{plane}
If $q\in X$ is a point of vanishing mean oscillation for  $\rho$,
then every tangent cone of $X$ at $q$ is an $n$-plane.
\end{proposition}

Without loss of generality, we
make a translation and  assume that $q=0$ in Proposition \ref{plane}.
Then fix a sequence $ r_i\to0$ such that \begin{equation}\label{xinfty}
X_i:=X_{0, r_i}\to X_{\infty}
\end{equation}
locally in the Hausdorff distance, where $X_\infty$ is a closed subset of $\R^N$. Our task is to show that
$X_\infty$ is an $n$-dimensional plane.

Write \begin{equation}\label{rhoi}
\rho^i(x):=\rho(r_ix)
\end{equation}
for each $i\in \N$. Then each $\rho^i$ is an $n$-tuple of flat $1$-forms defined
in $r_i^{-1}O$ with restrictions on $X_i$ as explained in \ref{intrinsic}. As before, we continue to denote these restrictions by $\rho^i$.

\begin{lemma}\label{taso1}
For every $R>0$ there exists  $\alpha \in (\wedge^1\R^N)^n$,
$\wedge\alpha\ne0$,
such that, upon passing to a subsequence
of $(r_i)$, we have \begin{equation}\label{inmeasure2}
\cH_{n}(\{x\in B(0,R) \cap X_{i}:||\rho^{i}(x)-\alpha||>2^{-i}\}) <
2^{-i}
\end{equation}
for every $i\in \N$.
\end{lemma}

\begin{proof}
Fix $R>0$. By the VMO-condition \eqref{vmo2}, we have that
\begin{equation}\label{vmo2}
\barint_{B(0,R) \cap X_i}||\rho^i-\rho_{B(0,r_iR) \cap X}|| \, d \cH_n =
\barint_{B(0,r_iR) \cap X}||\rho-\rho_{B(0,r_iR) \cap X}|| \, d \cH_n \to 0
\end{equation}
as $i\to\infty$. By passing to a subsequence
 of $(r_i)$, we may assume that \begin{equation}\label{vmo3}
\rho_{B(0,r_iR) \cap
X}=\barint_{B(0,r_{i}R) \cap X}\rho \, d \cH_n \to \alpha\,
\end{equation}
as $i\to\infty$, where $\alpha\in  (\wedge^1\R^N)^n$.
 By \eqref{cw2} and \eqref{cw3}, we have
  moreover that
  \begin{equation}\label{abound}
  0\,<\,c\,\le ||\wedge \alpha||\,\le C\,<\,\infty\,
  \end{equation}
  with constants depending only on local data.
  Now \eqref{vmo2} and \eqref{vmo3} jointly with uniform (local) Ahlfors regularity of
  $X_i$ give that \begin{equation}\label{vmo4}
\int_{B(0,R)\cap X_{{i}}}||\rho^{i}-\alpha|| \, d \cH_n\,\to\,0\,
\end{equation}
 as $i\to\infty$. By passing to a further subsequence, we obtain (\ref{inmeasure2})
 from \eqref{vmo4}, as required.
 \end{proof}

\begin{proof}[Proof of Proposition \ref{plane}]
Fix $R>0$. We use the subsequence and $\alpha$ given in Lemma \ref{taso1}.
Recall that $X$ is assumed   metrically oriented by a function $\xi$ as in \eqref{orient}.  Under the scaling $X\mapsto X_i$, we get  corresponding
orientations $\xi^i$ on $X_i$.
 Let $\zeta^i$ be the unique representative  of
the orientation $\xi^i$ on $X_i$ that corresponds to $\wedge\rho^i$ under the (Euclidean) duality (see  \ref{norms}).
Then we infer from (\ref{cw2}) and (\ref{cw3}) that
\begin{equation}\label{obs}
c_i(x)\,\zeta^i(x)\,=\,\xi^i(x)
\end{equation}
for almost every $x\in X_i$,
for some  measurable function
$c_i:X_i\to[c,C]$, where $0<c\le C<\infty$ are independent of $i\in\N$.
Next,  let $\zeta\in\wedge_n\R^N$ correspond to $\wedge\alpha\in\wedge^n\R^N$ under the duality. It follows that
\begin{align*}
\int_{B(0,R)\cap X_{{i}}}||\zeta^{i}-\zeta|| \, d \cH_n\,&\le C(n,N)\,\int_{B(0,R)\cap X_{{i}}}||\wedge\rho^{i}-\wedge\alpha|| \, d \cH_n\\
&\le\,C(n,N)\,\int_{B(0,R)\cap X_{{i}}}||\rho^{i}-\alpha|| \, d \cH_n\,\to\,0\,,
\end{align*}
where \eqref{vmo4} was used in the last step.
In particular, we find that
 \begin{equation}\label{vmo5}
 \int_{B(0,R)\cap X_{{i}}}||\xi^{i}-c_i\zeta|| \, d \cH_n\,\to\,0\,.
\end{equation}

By Proposition \ref{prop:locr}, we can  pass to another subsequence and find that $(X_i)$ converges in the sense of currents to
some locally rectifiable $n$-dimensional current $T$ with $\partial T = 0$. We  claim that the  support of the restriction $T\lfloor B(0,R)$ is equal to
$B(0,R)\cap X_\infty$ in the following sense:
as a rectifiable $n$-current $T\lfloor B(0,R)$ is represented by the rectifiable set
$B(0,R)\cap X_\infty$ together with a measurable orientation by a unit tangent $n$-vector field
$\eta: B(0,R)\cap X_\infty\to \wedge_n\R^N$ and an integrable almost everywhere nonzero
integer valued function $\theta: B(0,R)\cap X_\infty\to \bf Z$ \cite[4.1.28]{fed:gmt}.

To prove this claim, we first note that the support of $T$ is obviously a subset of $X_\infty$. Therefore, we only need to prove that the multiplicity function $\theta$ is almost everywhere different from zero.
To this end,
let $\psi: \R^N \longrightarrow \R$ be
a smooth nonnegative  function supported in $B(0,R)$.
We obtain from \eqref{vmo5} that
\begin{equation*} \begin{split}
 &\langle T,\psi\cdot\wedge \alpha\rangle=\lim_{i \to \infty} \langle X_i, \psi\cdot\wedge\alpha\rangle =
\lim_{i \to \infty} \int_{B(0,R)\cap X_i} \langle\psi\cdot\wedge\alpha,\xi^i \rangle d \cH_n\\
&\ge\liminf_{i\to\infty}\int_{B(0,R)\cap X_i} \langle\psi\cdot\wedge\alpha,\xi^i -c_i\zeta\rangle d \cH_n+ \liminf_{i\to\infty}\int_{B(0,R)\cap X_i} \psi\cdot\langle\wedge\alpha,c_i\zeta \rangle d \cH_n\\
&= \liminf_{i\to\infty}\,\,\int_{B(0,R)\cap X_i} \psi \cdot c_i\, d \cH_n\,\ge\,c\, \liminf_{i\to\infty}\,\int_{B(0,R)\cap X_i} \psi\, d \cH_n\,\\
&=\,c\,\int_{B(0,R)\cap X_\infty} \psi\, d \cH_n\,\ge\,c\,\cH_n(\{x\in B(0,R)\cap X_\infty:\psi(x)\ge1\})\,.\\
\end{split} \end{equation*}
Here in the last equality we also used the fact that the measures $\cH_n\lfloor X_i$ converge weakly
to $\cH_n\lfloor X_\infty$ (upon passing to another subsequence if necessary; see \cite[8.6]{ds:fractals}). If $\theta=0$ in a compact set $F$ of positive $\cH_n$ measure in $B(0,R)\cap X_\infty$, then we can choose a decreasing sequence of smooth functions $\psi_j$ with compact support
in $B(0,R)$ such that $\psi_j$ converges pointwise to the characteristic function of $F$. We then deduce from the preceding that
$$
0=\lim_{j\to \infty}\langle T,\psi_j\cdot\wedge \alpha\rangle\,\ge\,c\,\cH_n(F)>0\,,
$$
which is absurd. This proves the claim.

An argument similar to that in the preceding paragraph shows that
 $T\lfloor B(0,R)$ has  vanishing action
on forms with range orthogonal to $\wedge\alpha$.
This implies that almost every
approximate tangent plane of the set associated with
 $T\lfloor B(0,R)$, and therefore almost every
tangent plane of $B(0,R)\cap X_\infty$, is the $n$-plane
$V_\zeta$ determined by $\zeta$.

 To see that
$ B(0,R)\cap X_\infty$ is contained in a single plane, we use the Poincar\'e inequality  (\ref{poin}) valid in the limit space $X_\infty$  (see \cite[Section 9]{ch:lipschitz},
\cite[Section 8]{keith:zeit}, \cite{kos:pisa}). Fix $x_0\in B(0,\lambda^{-1}R)\cap X_\infty$, where $\lambda\ge1$ is as in \eqref{poin},
such that $T_{x_0}X_\infty$ exists and is equal to $V_\zeta$. Define
 $u: \R^N \to \R$ by
$$
u(x)= \dist( x, V_\zeta + x_0)\,.
$$
Then  $u$ is Lipschitz with vanishing
approximate differential almost everywhere along
$B(0,R)\cap X_\infty$. From the Poincar\'e inequality we obtain that $u$ is constant on $B(0,\lambda^{-1}R)\cap X_\infty$.
Since $u(x_0)=0$ we conclude that $u\equiv 0$ in $B(0,\lambda^{-1}R)\cap X_\infty$.
Because $R>0$ was arbitrary, we obtain that
$X_\infty$ is contained in the $n$-plane $V_\zeta$.

Finally,  because $\partial T=0$ and because the support of $T$ is contained in $V_\zeta$, we have in fact that the support of $T$
coincides with $V_\zeta$ (the constancy theorem \cite[p.\ 357]{fed:gmt}). In particular, $X_\infty=V_\zeta$, and the
 proof of  Proposition \ref{plane} is thereby complete.
\end{proof}

\begin{remark}\label{armk}
The proof of Proposition \ref{plane} shows that each of the limit  $n$-tuples of 1-covectors $\alpha\in (\wedge^1\R^N)^n$,
corresponding to a fixed value $R>0$ and to a choice of a subsequence as in the proof of Lemma \ref{taso1}, determines an orientation on the limit $n$-plane $X_\infty$ by the dual of $\wedge \alpha$. We call such an orientation on a
tangent $n$-plane at a VMO-point
a {\it limiting orientation}.

Note that {\it a priori} there could be two such limiting orientations
on a single $X_\infty$.
It follows from Lemma \ref{Lisbld} below that after passing, initially, to a certain subsequence of $(X_i)$, such an orientation is in fact independent of $R$.
\end{remark}

For the next lemma, assume that $q\in X$ is
a point of vanishing mean oscillation for $\rho$ and that $X_i=X_{q,r_i}$ converges to a tangent $n$-plane $X_\infty$ such that \eqref{inmeasure2} holds for some $R$ and $\alpha$. Assume further that $X_\infty$ is equipped with a limiting orientation corresponding to $\alpha$
as explained in Remark \ref{armk}. Let $\pi:\R^N\to X_\infty$ be the orthogonal projection, and abbreviate $\pi_i:=\pi|X_i$. The local degree $\mu(0,\pi_i, D)\in \bf Z$
is defined for every $i$ and for every connected open set $D\subset X_i$ with compact closure in $X_i$ such that
$0\notin \pi_i(\partial D)$.  (See \cite[Chapter I.4]{rick:quasiregular}
and \cite[Section 2]{hr:duke} for the definition of the local degree for  mappings between oriented homology manifolds.) For $s>0$,
denote by $D_i(s)$
the component of $B(0,s)\cap X_i\subset X_i$ that contains $0\in X_i\cap X_\infty$. For every fixed $s$, the local degree
$\mu(0, \pi_i,D_i(s))$ is defined for all sufficiently large $i$.

\begin{lemma}\label{deglemma0}
For every $0<s<R/2$  the local degree  $\mu(0, \pi_i, D_i(s))$ is positive
for all sufficiently large $i$.
\end{lemma}

\begin{proof}
After some preparatory work, the proof is similar to that in \cite[Proposition 6.17]{hr:duke}.
First, denote by $G_i(s)$ the set of those points $x\in D_i(s)$ such that $T_xX_i$ exists and that the differential $D\pi_i(x):T_xX_i\to X_\infty$ has rank $n$ and is sensepreserving with the given orientations. It follows from \eqref{inmeasure2} that
\begin{equation}\label{degk1}
\cH_n(D_i(s)\setminus G_i(s))\to 0\,, \ \ \ \ \ \ \ i\to \infty\,.
\end{equation}
Second, the change of variables formula (see \cite[3.2.20]{fed:gmt} or \cite[(5.5)]{hr:duke}) implies that
\begin{equation}\label{degk2}
\int_{\pi_i(D_i(s))}N(\pi_i, D_i(s), y)\,d\cH_n(y)\,\le\, C\,s^n\,,
\end{equation}
where $N(\pi_i, D_i(s), y):={\text{card}}\{\pi_i^{-1}(y)\cap D_i(s)\}$ and $C>0$ depends only on the local data near $q$. Third, we claim that \begin{equation}\label{degk3}
B(0,s/C)\cap X_\infty\subset \pi_i(D_i(s))
\end{equation}
for all large $i$, where $C\ge 1$ depends only on the local data. To prove \eqref{degk3}, we need to observe that there is $0<\lambda <1$, depending only on the local data, such that \begin{equation}\label{degk4}
B(0,\lambda s)\cap X_i\subset D_i(s)
\end{equation}
for all large $i$; this follows from the local linear contractibility. Then we
use the maps
$\psi_{\infty,i}:B(0,R)\cap X_\infty\to B(0,2R)\cap X_i$ as in \eqref{retrk1}. In fact, for all large $i$, we have that $\pi_i\circ\psi_{\infty,i}:B(0,2R)\cap X_\infty\to X_\infty$ is homotopic to the identity via a radial homotopy $H_i^t$ such that $0\notin H_i^t(\partial B(0,\lambda s)\cap X_\infty)$, where $\lambda$ is as in \eqref{degk4}. In particular, it follows from the properties of the local degree that the image of  $\pi_i\circ\psi_{\infty,i}$, and hence the image of $\pi_i$, has to contain $B(0,\lambda s/2)\cap X_\infty$.
Because $X_i\to X_\infty$, we deduce from \eqref{degk4} that indeed $B(0,\lambda s/2)\cap X_\infty$ is contained in $\pi_i(D_i(s))$. Thus \eqref{degk3} follows.

By combining \eqref{degk1} -- \eqref{degk4}, we infer that for every sufficiently large $i$ there is a point $y_i$ in the 0-component of $X_\infty\setminus \pi_i(\partial D_i(s))$ such that $\pi_i^{-1}(y_i)\cap D_i(s)$ consists of finitely many points only, all belonging to $G_i(s)$. The argument
now runs similarly to that in \cite[Proposition 6.17]{hr:duke}, and we conclude the proof of the lemma.
\end{proof}

We have the following corollary, valid in the situation of Lemma \ref{deglemma0}.

\begin{corollary}\label{deglemma0koro}
The local degree $\mu(0,\psi_{\infty,i}, B(0,s)\cap X_\infty)$ is positive for all sufficiently large $i$, where $\psi_{\infty,i}$ are as in \eqref{retrk1}.
\end{corollary}

\section{Tangent maps}\label{tmaps}
We assume the basic set-up as in \ref{setup}.

Given $E\subset \R^N$, $x \in E$, and a function $g:E \to \R^n$,
define $$
g_{x,r}: E_{x,r}
\to \R^n\,
$$ by
$$
g_{x,r} (y) := \frac{g(x+ry)- g(x)}{r}\,,
$$ where we recall $E_{x,r}=r^{-1}(E-x)$.
Further, given a sequence of sets $(E_i)$ together with functions $g_i:E_i \to \R^n$,
we say that $(g_i)$ \emph{converges to} $g$, where $g:E\to \rn$, if
$(E_i)$ converges to $E$ locally in the Hausdorff distance
 and if for every sequence $(y_i)$ of points
converging to $y $, where $y_i \in E_i$ and $y\in E$, we have that $g_i(y_i)$
converges to  $g(y)$. We say that $(g_i)$ \emph{converges to} $g$ \emph{uniformly on bounded sets}
if, in addition, the convergence $g_i(y_i)\to g(y)$ is uniform whenever the points $y_i, y$ lie in a bounded set in $\R^N$. (Compare \cite[8.3]{ds:fractals}.)

Recall from \eqref{omap} the notation $F_q(x)$
for $q\in X$ and $x\in O$.

\begin{proposition} \label{linea}
Let $q\in X$ be a  point of vanishing mean oscillation for $\rho$ and let $r_i\to0$ be a sequence such that $X_{q,r_i}\to X_\infty$, where {\rm{(}}necessarily, by Proposition \ref{plane}{\rm{)}} $X_\infty$ is an
$n$-plane.
Then  there is
 a subsequence $(r_{i_j})$ of $(r_i)$ such that
 the maps
$$
(F_q)_{q,r_{i_j}}:X_{q,r_{i_j}}\to \rn
$$
converge uniformly on bounded sets to
a bijective linear map $L:X_\infty\to \rn$ as $r_{i_j}\to\infty$. Moreover, $L$ is bi-Lipschitz
with constant depending only on the local data.
\end{proposition}

Proposition \ref{linea} will be proved with the assistance of
several lemmas. In the course of the proof, we will pass to various subsequences
without renamings.

Without loss of generality, we
make a translation and assume that $q=0$.
Write $X_i:=X_{0,r_i}$
and \begin{equation}\label{Fi}
F^i: = (F_0)_{0,r_i}\,, \ \ \ \ \ \ F^i:X_i\to \R^n\,.
\end{equation}
Then $(F^i)$ is a collection of BLD-mappings with uniform data, as explained in subsection \ref{cartBLD}. In particular, each $F^i$ is uniformly
Lipschitz and by an Arzel\`a-Ascoli type argument (see \cite[Lemma 8.6]{ds:fractals})
we can  pass to a subsequence
so as to guarantee that $(F^i)$ converges uniformly on bounded sets to some Lipschitz map $L:X_\infty\to \rn$. By passing to a yet another subsequence, we may assume that $X_\infty$ is equipped with a limiting orientation as in Remark \ref{armk}.

Finally, note that \eqref{omap} defines $F^i(x)$ for every $x\in B(0,R)$ for a fixed $R>0$ and $i$ large enough. Then $F^i:B(0,R)\to \R^n$ is Lipschitz with constant independent of $i$ \cite[4.6]{hs:duke}, and we may assume that $F_i$ converges uniformly on bounded sets to a Lipschitz map $\R^N\to \R^n$. We call this map $L$ as well.

\begin{lemma}\label{Lisbld}
The mapping $L:X_\infty\to \rn$ is an orientation preserving BLD-mapping,
where  $X_\infty$ is equipped with a limiting orientation.
\end{lemma}

\begin{proof}
Note first that $L$ is Lipschitz. It follows from \cite[Theorems
4.5 and 6.8]{hr:duke}
 that in the present context BLD-mappings can be characterized as
sensepreserving {\it regular maps}
 (as in \cite[Definition 12.1]{ds:fractals}). Therefore, the proposition
 follows from the stability of regular
maps under limits \cite[Lemma 12.7]{ds:fractals}, provided $L$ is sensepreserving. This will be shown below.
(Note that in the preceding reasoning, we need a routine localization of
\cite[Theorem 4.5]{hr:duke} and also need to know that the spaces $X_i$ are appropriately quasiconvex.
The latter fact follows from the validity of a uniform Poincar\'e inequality on $X_i$  as explained in subsection \ref{sobolev} and from \cite[Section 17]{ch:lipschitz}.)

Next, we invoke the following uniform estimate for the local distortion of BLD-mappings \cite[Proposition 4.13]{hr:duke}: for every $R>0$ there exists $i_0$ such that \begin{equation}\label{locdist}
B(F^i(x), r/C)\subset F^i(B(x,r)\cap X_i)\subset B(F^i(x), Cr)
\end{equation}
whenever $i\ge i_0$, $x\in B(0,R)\cap X_i$, and $0<r<R$. (That the assumptions in \cite[Proposition 4.13]{hr:duke} are satisfied in the present situation is explained as in the preceding paragraph.)
It easily follows from \eqref{locdist} that $L$ is an open mapping.
Discrete and open mappings between manifolds are either sensepreserving or sensereversing \cite[5.2]{vai:minimal}. It suffices, therefore, to show that
the local degree $\mu(0,L, B(0,s)\cap X_\infty)$ is positive, where the radius $s>0$ is so small that
$f^{-1}(0)\cap B(0,2s)=\{0\}$.
(See \cite[Chapter I.4]{rick:quasiregular}
and \cite[Section 2]{hr:duke} for the definition of the local degree.)

To this end, fix $s>0$ as above. We use the maps $\psi_{\infty,i}:X_\infty\cap B(0,R)\to X_i$ as in Proposition \ref{retr} for $R>s>0$ fixed and for $i$ sufficiently large.
It follows from \eqref{locdist} that the maps $F^i\circ\psi_{\infty,i}|B(0,s)\cap X_\infty$ and $L|B(0,s)\cap X_\infty$ are homotopic through maps
$H_i^t$ such that $0\notin H_i^t(\partial B(0,s)\cap X_\infty)$ for all  $0\le t\le1$, for all
large $i$. In particular, the local degrees satisfy
$$
\mu(0,F^i\circ\psi_{\infty,i},B(0,s)\cap X_\infty)=\mu(0,L,B(0,s)\cap X_\infty)
$$
for all large $i$. Since $F^i$ is sensepreserving and $\mu(0, \psi_{\infty,i}, B(0,s)\cap X_\infty)$ is positive by Corollary \ref{deglemma0koro}, we have that
$\mu(0,L,B(0,s)\cap X_\infty)$ is positive as required.
This completes the proof of the lemma.
\end{proof}

For the remainder of this section, we fix $R>0$ and pass to a  subsequence as in Lemma \ref{taso1}. In particular, we assume that
\eqref{inmeasure2} holds for  an $n$-tuple of covectors $\alpha \in (\wedge^1\R^N)^n$ with $\wedge\alpha\ne0$.

We require the concept of the $n$-modulus of a curve
family in metric measure spaces \cite{hk:quasi},  \cite{hei:lectures}.

\begin{lemma}\label{loewner}
 There exists a constant $C_1\ge1$  depending only on
the local data such that following holds for all large $i$: every collection
$\Gamma$ of curves  with length at most $C_1 |x-y|$ connecting $B(x, r)\cap X_i$
to $B(y,r)\cap X_i$ in $X_i$, where $x,y \in B(0, R)\cap X_i$ and $0<2r<|x-y|$, satisfies
$\bmod_n (\Gamma) \ge C_2$, where $C_2>0$ depends only on the local data and on
the ration $|x-y|/r$.
\end{lemma}

\begin{proof} The validity of a (local) Poincar\'e inequality on $X$ (see \cite{sem:curves} and \ref{sobolev}) implies that $X$ is (locally) a {\it Loewner space}.
 Thus, the claim follows from \cite[Lemmas
3.15 and 3.17]{hk:quasi}.
\end{proof}

Now fix $z\in X_\infty\cap B(0,R/2C_1)$, where $C_1$ is as in Lemma \ref{loewner}.
We only consider indices $i$ large enough such that the conclusion of Lemma \ref{loewner} holds.

In the next technical lemma,
for each $i\ge1$, we subdivide the line segment
$[0,z]$ into $2^i$ line segments of equal length $2^{-i}|z|$
and denote by $z_{i,j}$, $j=0,1,\dots,
2^{i}-1$, the beginning points of these segments
naturally ordered with $j$ such that $z_{i,0}=0$. We also put
 $z_{i,2^i}=z$.

\begin{lemma} \label{flata}
After passing to a subsequence of $(r_i)$, there exists a sequence
of rectifiable curves $(\gamma_i)$ parametrized by the arc length,
$\gamma_i:[0,\leng(\gamma_i)]\to X_i$, such that $\gamma_i(0)\to 0$, that $\gamma_i(\leng(\gamma_i))\to z$, and that
each $\gamma_i$ decomposes into $2^i$ consecutive subarcs
$\gamma_{i,j}$, $j=0,1,\dots,
2^{i}-1$, with the following two properties:

{\rm (a)} denoting the beginning and the end points of $\gamma_{i,j}$ by
$a_{i,j}$ and $b_{i,j}$, respectively, where $$
0\le \gamma_{i,j}^{-1}(a_{i,j})<\gamma_{i,j}^{-1}(b_{i,j})=\gamma_{i,j+1}^{-1}(a_{i,j+1})
< \leng(\gamma_i)\,,
$$
we have that \begin{equation}\label{flatak}
|a_{i,j}-z_{i,j}|<2^{-i^2}|z|\,, \ \ \ \ |b_{i,j}-z_{i,j+1}|<2^{-i^2}|z|\,
\end{equation}
for each $j$;

{\rm (b)} the length of each $\gamma_{i,j}$ does not exceed $C\,2^{-i}|z|$,
where $C>0$ is a constant depending only on the local data.

In particular, we have that
$\leng(\gamma_i)\le C\,|z
|$, where $C>0$ depends
only on the local data.

Moreover, we have that
\begin{equation}\label{rtoal}
\lim_{i\to\infty}\int_{\gamma_i}|\rho^i-\alpha|\,ds\,=\,0\,,
\end{equation}
where  the line integral in  \eqref{rtoal} should be understood
via the canonical restrictions of forms to 1-rectfifiable sets as explained in \ref{rest}.
\end{lemma}

\begin{proof}
Fix $m\ge3$. Since $X_i\to X_\infty$, where $X_\infty$ is an
 $n$-plane, there exists $M_1\ge m$
such that $$
B(z_{m,j}, 2^{-m^2-1}|z|)\cap X_i \neq
\emptyset
$$
for each $j=0, 1, \dots, 2^m$ and each $i\ge M_1$. Fix such an $i$.
 By Lemma \ref{loewner},
the collection
$\Gamma_{i,j}$ of curves with length at most
$C_12^{-m}|z|$ connecting $ B(z_{m,j},2^{-m^2}|z|)\cap X_i$ to
 $B(z_{m,j+1},2^{-m^2}|z|)\cap X_i$
satisfies
\begin{equation}\label{suuri}
\bmod_n (\Gamma_{i,j})\ge C_2>0\,,
\end{equation}
for each $j=0,1,\dots, 2^m-1$,
where $C_2>0$ depends only on the local data and $m$.

We claim that there exists $M_2\ge \max\{M_1, m^2+1\}$ with the following property:
if $i\ge M_2$, then for each $j=0,1,\dots, 2^m-1$ one finds a curve
$\gamma_{i,j}\in \Gamma_{i,j}$ satisfying
\begin{equation}\label{vahan}
\cH_1(\gamma_{i,j}\cap\{x:||\rho^i(x)-\alpha||\ge 2^{-m^2}\})\,\le\,2^{-m^2}\,.
\end{equation}
Indeed, if $$
\cH_1(\gamma\cap\{x:||\rho^i(x)-\alpha||\ge 2^{-m^2}\})\,>\,2^{-m^2}\,
$$
for all $\gamma\in\Gamma_{i,j}$, then it follows from the definition of modulus and from \eqref{inmeasure2}
that
$$
\bmod_n(\Gamma_{ij})\,\le\,2^{2m^2n}\,
\mathcal H_n(\{x\in B(0,R)\cap X_i:||\rho^i(x)-\alpha||\ge 2^{-m^2}\})\,<\,2^{2m^2n}\,2^{-i}\,.
$$
This contradicts \eqref{suuri}, provided $i$ is large enough. (Note that here we used the fact
that the curves in
each $\Gamma_{i,j}$ lie in $B(0,R)$, which is true  by the choices.)

Now assume that $i\ge M_2$ and pick a curve $\gamma_{i,j} \in \Gamma_{i,j}$ such that \eqref{vahan} holds. By
 the quasiconvexity of $X_i$ (see \cite[Section 17]{ch:lipschitz}), we can extend $\gamma_{i,j}$ by a curve of length at most $C\,2^{-m^2}|z|$
so that
the end point of $\gamma_{i,j}$ coincides with the
beginning point of $\gamma_{i,j+1}$ (keeping
 with the same notation for the new extended curves).
Here $C>0$ depends only on the local data. Finally (see
\cite[Section 2]{vai:exh}, for example)
we can assume that the curves $\gamma_{ij}$ are arcs, or embeddings
as maps by using the arc length parametrizations.
For these arcs we have that
\begin{equation}\label{vahan2}
\cH_1(\gamma_{i,j}\cap\{x:||\rho^i(x)-\alpha||\ge 2^{-m^2}\})\,\le\,2^{-m^2}+C\,2^{-m^2}|z|\,,
\end{equation}
where $C>0$ depends only on the local data.

We let $\gamma_i$ be the consecutive union of
these  $\gamma_{i,j}$'s over $j$. Upon obvious renaming, we thus have a sequence that
 satisfies the
required properties (a) and (b) by construction, while  \eqref{rtoal}
follows from \eqref{vahan2}.
The lemma is thereby  proved.
\end{proof}

We  view the curves $\gamma_i$
in the preceding lemma as one dimensional rectifiable currents in a canonical way,
\begin{equation}\label{cancur}
\gamma_i=\sum_{j=0}^{2^i-1}\gamma_{i,j}\,,
\end{equation}
where each of the arcs $\gamma_{i,j}$ are interpreted as rectifiable currents
with unit multiplicity and orientation coming from $\gamma_i$ traversing
 from $a_i:=\gamma_i(0)$ to $b_i:=\gamma_i(\leng(\gamma_i))$.
In particular, we can view the curves $\gamma_i$ as flat 1-chains.

\begin{lemma}\label{flgconv}
The sequence of currents
$(\gamma_i)$ converges to
$[0,z]$ in the flat norm.
\end{lemma}

\begin{proof}
By using the decomposition of each $\gamma_i$ guaranteed by Lemma \ref{flata},
we write
$$
[0,z]-\gamma_i=[0,a_i]-[z,b_i]+\sum_{j=0}^{2^i-1}T_{i,j}\,,
$$
where
$$
T_{i,j}:=[z_{i,j},z_{i,j+1}]+[z_{i,j+1},b_{i,j}]-\gamma_{i,j}-[z_{i,j},a_{i,j}]\,
$$
is a rectifiable closed current. By the isoperimetric
inequality \cite[4.2.10 and 4.5.14]{fed:gmt}  there is a flat 2-chain $S_{i,j}$ such that $\partial S_{i,j}=T_{i,j}$
and that
$$
4\pi|S_{i,j}|\,\le\,|T_{i,j}|^{2}\,,
$$
where $|\cdot|$ stands for mass.
Therefore,
\begin{align*}
|[0,z]-\gamma_i|_\flat&\,\le\,|[0,z]-\gamma_i-\partial
\sum_{j=0}^{2^i-1}S_{i,j}|+|\sum_{j=0}^{2^i-1}S_{i,j}|\\
&\le\,|0-a_i|+|z-b_i|+\sum_{j=0}^{2^i-1}|T_{i,j}|^2\\
&\le\,2^{-i^2}|z|+2^{-i^2}|z|+C\,2^{-i}|z|^2\,,\\
\end{align*}
where $C>0$ depends only on the local data.
This proves the lemma.
\end{proof}

We are now ready for the final lemma.

\begin{lemma}\label{Linear}
We have that
\begin{equation}\label{Lq}
L(z) = \int_{[0,z]} \alpha\,
\end{equation}
for every $z\in B(0, R/2C_1)\cap X_\infty$, where $C_1>0$ is as in Lemma \ref{loewner}.
In particular, we have that $\alpha$ is independent of $R$.
\end{lemma}

\begin{proof}
Recall that $R>0$ is fixed and that we wish to prove \eqref{Lq} for $z\in B(0,R/2C_1)$, where $L$ is the limit of the maps $F^i$ as in \eqref{Fi}. Also recall that $F^i$ is defined at $z$ for $i$ large enough, as explained before Lemma \ref{Lisbld}.
Note the several passages to subsequences in the preceding lemmas.
 We use the curves $\gamma_i$ as guaranteed by Lemma \ref{flata}.

The preceding understood, we estimate
$$ |F^i(z)-\int_{\gamma_i}\rho^i| =
|\int_{[0,z]-\gamma_i} \rho^i| \,\le\,||\rho^i||_\flat\,
|[0,z]-\gamma_i|_\flat\,.
$$
Because the flat norms of the forms $\rho^i$ are uniformly bounded,
we therefore obtain from Lemma \ref{flgconv} that
$$
\lim_{i\to\infty}\int_{\gamma_i}\rho^i=\lim_{i\to\infty}F^i(z)=L(z)\,.
$$
On the other hand, \eqref{rtoal} and Lemma \ref{flgconv}
give that
$$
\lim_{i\to\infty}\int_{\gamma_i}\rho^i=
\lim_{i\to\infty}\int_{\gamma_i}\alpha\,=\,\int_{[0,z]}\alpha\,.
$$
In conclusion, we obtain \eqref{Lq}, as required.
\end{proof}

Proposition \ref{linea} obviously follows from Lemma \ref{Linear} and from \eqref{abound}.

\section{Proof of Theorem \ref{main2}}\label{pofmain}
In this section, we complete the proof of our main result Theorem \ref{main2}.
We assume the notation and conventions made in \ref{setup}. In addition, {\it we assume now that
$\rho\in H^{1,2}(X)$.}

We consider  maps
$F_q:X\to \R^n$ as in \eqref{omap},  for $q \in X$. Recall however that $F_q(x)$ is defined for all $x\in O$, where $O$ is a fixed convex $\R^N$-neighborhood of $p$.
 By the remarks made in Section \ref{cartBLD}, to prove Theorem
\ref{main2} it suffices to show that $p$ does not belong to the branch set $B_{F_p}$ of $F_p$. We suppose towards a contradiction that $p \in B_{F_p}$.
To achieve a contradiction, we will analyze limits of rescaled maps and spaces
near $p$ by using the results in the previous section.

As earlier, we require several lemmas.

\begin{lemma}\label{n-2}
There  exists  a sequence $(q_k)$ of points
in the branch set $B_{F_p}$ converging to $p$ such that each point $q_k$ is a point of vanishing mean oscillation  for $\rho$.
\end{lemma}

\begin{proof}
It follows from the Sobolev assumption  $\rho\in H^{1,2}(X)$ and from the Poincar\'e inequality \eqref{poin} that
\begin{equation}\label{poin2}
\left(\barint_{B \cap X }|\rho-\rho_B|\, d \cH_n\right)^2 \le \,C\,
({\text{diam}}B)^2\barint_{ \lambda B \cap X}|D\rho|^2 \, d \cH_n
\end{equation}
for each small $\R^N$-ball $B$ near $p$, where $\lambda, C\ge1$ are
independent of $B$.
By using \eqref{poin2} and  standard covering
theorems as in  \cite[Chapter 1]{hei:lectures} or \cite[Theorem 3, p. 77]{eg:measure}, for example, we obtain that $\cH_{n-2}$-almost every point is a point of vanishing mean oscillation for $\rho$.
On the other hand, as discussed in Section \ref{cartBLD}, we have that
$\cH_{n-2}(B_{F_p}\cap V)>0$
for every neighborhood $V$ of $p$.
The lemma follows from these remarks.
\end{proof}

To continue the proof of Theorem \ref{main2}, fix a sequence $(q_k)$ as in Lemma \ref{n-2}.
We may assume that $B(q_k, |p-q_k|)\subset O$ for all $k$.
By Proposition \ref{linea}, we can find for each $k\ge1$  a
sequence $(r_{k,i})$ of positive reals, converging to naught as $i\to\infty$, such that the maps
$$
(F_{q_k})_{q_k,r_{k,i}}:X_{q_k,r_{k,i}}\to\R^n
$$
converge, as $i \to
\infty$, uniformly on bounded sets to a bijective linear map $L_k:P_k\to\rn$,  where $P_k$ is an $n$-plane in $\R^N$ with $X_{q_k, r_{k,i}}\to P_k$.
 Each of the $n$-planes $P_k$ carry a natural limiting orientation inherited from $X$, and
the linear maps $L_k:P_k\to \R^n$ are orientation preserving and bi-Lipschitz with constant independent of $k$.
(See   Remark \ref{armk} and Lemma \ref{Lisbld}.)
We also assume that each $L_k$ is defined in all of $\R^N$ and is Lipschitz with uniform Lipschitz constant.

From now on, we understand the sequences $(r_{k,i})$ as fixed.

\begin{lemma}\label{lah}
For every $R>0$ and every $x\in X_{q_k,r_{k,i}}$ with $|x|\le R$ we have that \begin{equation}\label{lop1}
| (F_{q_k})_{q_k,r_{k,i}} (x)-(F_p)_{q_k,r_{k,i}} (x) |\,\le\,\,||d\rho||_\infty \,|x|\,|p-q_k|
\end{equation}
whenever $r_{k,i}R<|p-q_k|$.
\end{lemma}

\begin{proof}
Note that it follows from the choices that both maps in \eqref{lop1} are defined for every $R$ and $x$ as required. Thus, fix $R>0$ and  $x\in X_{q_k,r_{k,i}}$ such that $|x|\le R$, and assume that $r_{k,i}R<|p-q_k|$.
 Then
\begin{align*}
| (F_{q_k})_{q_k,r_{k,i}} (x)-(F_p)_{q_k,r_{k,i}} (x) |
&= \frac{1}{r_{k,i}} \left|\int_{[q_k,q_k+r_{k,i}x]}\rho-
\int_{[p,q_k+r_{k,i}x]}\rho+\int_{[p,q_k ]} \rho \right|\\
& = \frac{1}{r_{k,i}} \left| \int_{[p, q_k, q_k+r_{k,i}x]} d \rho \right| \\
&\le \frac{1}{r_{k,i}} \|d\rho\|_\infty\, |r_{k,i}x| \, |p -
q_k|\,=\,||d\rho||_\infty\,|x|\,|p-q_k|\,,
\end{align*}
where $[p,q_k,q_k+r_{k,i}x]$ denotes the (oriented) triangle span by the three points.
The lemma follows.
\end{proof}

Next, fix linear isometries
$\tau_k:\R^N\to \R^N$ such that $\tau_k|\rn:\rn\to P_k$
is an orientation preserving surjection. Here and later we view $\R^n$ as a subset of $ \R^N$ in a natural way.    The maps $L_k$ are uniformly bi-Lipschitz,
so after passing to a subsequence
we can further assume that the maps $L_k\circ\tau_k|\R^n:\rn\to\rn$ converge to an invertible sensepreserving linear map $L:\rn\to\rn$ locally uniformly. Then the maps $(F_p)_{q_k,r_{k,i}}\circ\tau_k$
 are defined for every $x\in B(0,R)\cap \tau_k^{-1}(X_{q_k,r_{k,i}})$ such that $r_{k,i}R<|p-q_k|$.

\begin{lemma}
After passing to a subsequence of the sequence $(q_k)$, there
exists a sequence of positive real numbers $(s_k)$ converging to
zero such that  the sets $$
Y_k:=B(0,k)\cap \tau_k^{-1}(X_{q_k,s_k})
$$ converge locally in the Hausdorff distance to
$\R^n$ and such that the maps
$$
(F_p)_{q_k,s_k}\circ\tau_k:Y_k\to\R^n
$$
are defined and converge locally uniformly to $L$ as $k\to\infty$.
\end{lemma}

\begin{proof}
We may assume that the numbers $r_{k,i}$ satisfy $r_{k,i}k<|p-q_k|$ for all $k$ and $i$, so in particular
the maps $(F_p)_{q_k,r_{k,i}}\circ\tau_k$ are always defined on $B(0,k)\cap \tau_k^{-1}(X_{q_k,r_{k,i}})$.
Fix $m\ge1$. Because $L_k\circ\tau_k\to L$ locally uniformly, there exists $k_m\ge m$ such that
\begin{equation}\label{huh1}
|L(y)-L_k\circ\tau_k(y)|<\frac1{m}
\end{equation}
whenever $k\ge k_m$ and $y\in B(0,m)\cap \R^n$. Moreover,
 for every $k$ there exist $i_{k,m}$ and $\delta_{m}>0$ such that \begin{equation}\label{huh2}
|L_k\circ\tau_k(y)-(F_{q_k})_{q_k,r_{k,i}}\circ\tau_k(y_{k,i})|<\frac1{m}
\end{equation}
whenever $i\ge i_{k,m}$ and \begin{equation}\label{huh25}
y\in B(0,m)\cap \rn\,, \ \ \ \ \ y_{k,i}\in B(0,m)\cap  \tau_k^{-1}(X_{q_k,r_{k,i}})\,
\end{equation}
satisfy $|y-y_{k,i}|<\delta_{m}$. We can also assume, by Lemma \ref{lah}, that $k_m$ and
$i_{m,k}$ are large enough so that in addition to \eqref{huh1} and \eqref{huh2} we have \begin{equation}\label{huh3}
| (F_{q_k})_{q_k,r_{k,i}} \circ\tau_k(y_{k,i})-(F_p)_{q_k,r_{k,i}} \circ\tau_k(y_{k,i}) |\,<\,\frac1{m}\,
\end{equation}
whenever $k\ge k_m$, $i\ge i_{k,m}$, and $y_{k,i}$ is as in \eqref{huh25}.
Now consider  the subsequence $(q_{k_m})$ of $(q_k)$, and  the concomitant
 numbers
 $s_m:=r_{k_m, i_{k_m,m}}$ converging to zero. By passing to a yet new subsequence
 of $({k_m})$ and enlarging $i_{k,m}$, if necessary, we may assume that the sets
 $B(0,m)\cap \tau_{k_m}^{-1} (X_{q_{k_m}}, s_m)$ converge as required to $\R^n$. In consequence, it
 follows from \eqref{huh1} -- \eqref{huh3} that (upon obvious renaming) we have
 sequences
 as asserted. The lemma follows.
\end{proof}

We  use  the subsequence provided by the preceding lemma, and
set
\begin{equation}\label{fkoo}
f_k := L^{-1} \circ (F_p)_{q_k,s_k}\circ \tau_k\,.
\end{equation}
Thus, $f_k:Y_k\to \rn$
is a BLD-mapping with $f_k(0)=0$.  By construction, the sets $Y_k$ converge locally in the Hausdorff distance to $\rn$, and the mappings $f_k$ converge locally uniformly to the identity on $\rn$. For $s>0$, denote by $D_k(s)$ the 0-component of $B(0,s)\cap Y_k\subset Y_k$. Then the local degree $\mu(0, f_k, D_k(s))$ is defined for every sufficiently large $k$.
(Compare the discussion before  Lemma \ref{deglemma0}.)

\begin{lemma}\label{deglemma}
Given $s>0$, the local degree satisfies
\begin{equation} \label{dega}
\mu(0, f_k, D_k(s))=1 \, \end{equation}
for all sufficiently large $k$. \end{lemma}

\begin{proof}
Fix $s>0$.
As explained in Proposition \ref{retr}, for all sufficiently large $k$ there exist  mappings
$\psi_k :  B(0,4s)\to Y_k$ such that
$$
| \psi_k(x)-x| \le C \dist(x, Y_k)
$$
for every $x\in B(0,4s)$, where $C\ge
1$ is independent of $k$.
By pushing  the obvious radial homotopies in $\R^N$  to $Y_k$ by using the maps $\psi_k$,
we easily infer that
$$
 \psi_k\circ f_k|B(0,2s)\cap Y_k:B(0,2s)\cap Y_k\to Y_k
$$
are all homotopic to the identity through homotopies $H^t_k$ such that $0\notin H_k^t(\partial D_k(s))$,
for all large $k$.
The preceding understood, it follows from Corollary \ref{deglemma0koro} and from the properties of the local degree that (\ref{dega}) holds. \end{proof}

\begin{proof}[Proof of Theorem \ref{main2}]
Recall that $q_k \in B_{F_p}$,
 where $q_k$ is given in Lemma \ref{n-2}. In particular, we have that
 for each $k$ the mapping  $f_k:Y_k\to\R^n$ as defined in \eqref{fkoo} is
BLD and has branching at $0$.
This contradicts, for large $k$, the local degree value
given in \eqref{dega}. The proof of Theorem \ref{main2} is thereby complete.
\end{proof}

\subsection{Remarks}\label{vmormks}
As discussed in the introduction, the requirement that a Cartan-Whitney presentation
$\rho$ belongs to $H^{1,2}_{loc}(X)$  in Theorems \ref{main}, \ref{main3}, and \ref{main2}
 is a kind of  second order requirement (for the coordinates).
The proof  here and in earlier sections  shows that this hypothesis could be changed
to a VMO-type condition. Indeed, the Sobolev condition was used only in order to have the
conclusion of Lemma \ref{n-2}. For the same conclusion, it would suffice
to assume that
\begin{equation}\label{uusi}
{\text{ $\mathcal H_{n-2}${\it -almost every point in $X$ is a {\rm{VMO}}-point for} $\rho$}}.
\end{equation}
(Recall the definition for VMO-points from \eqref{vmo}.)
In particular, all  the main results of this paper would remain
valid if the hypothesis that $\rho$ is in $H^{1,2}$ is changed to requirement  \eqref{uusi}.

\section{Smoothability}\label{sbilitysec}
In this section, we first formulate and prove the following general smoothability result.

\begin{theorem}\label{sbility2}
Let $X\subset \R^N$ be a locally Ahlfors $n$-regular, connected,  locally linearly contractible, and oriented homology $n$-manifold. Assume that there is a Cartan-Whitney presentation $\rho=(\rho_1,\dots, \rho_n)$
for $X$ that is defined in an $\R^N$-neighborhood $W$ of $X$.
  If  the restriction of $\rho$ to $X$ belongs to the local Sobolev space $H^{1,2}_{loc}(X)$,
then $X$ is a smoothable topological $n$-manifold.
\end{theorem}

The proof of Theorem \ref{sbility2} is a combination of our Theorem \ref{main2} and the argument from \cite[p.\ 336]{sul:bld}. For completeness, we present the details. We begin by describing a canonical metric on $X$ induced by a globally defined Cartan-Whitney presentation. Such metrics were employed in \cite{sul:bld}, and the crucial point is to produce local bi-Lipschitz homeomorphisms with respect to this metric that are near isometries. A theorem of Shikata \cite{shi2:smooth} (and its noncompact version
in \cite{gm:smooth}) then gives smoothability.

\subsection{Metrics induced by Cartan-Whitney presentations}\label{rhodconst}
Let $X\subset W\subset \R^N$ and $\rho$ be as
in
Theorem \ref{sbility2}, except that no Sobolev regularity is required for the Cartan-Whitney presentation $\rho$. Recall that for every flat 1-form $\omega$ in an open set $V$ in $\R^N$, and for every Lipschitz map $\gamma:[0,L]\to V$, the action
\begin{equation}\label{eka0}
\langle \omega(\gamma(t)),\gamma'(t)\rangle\,\le\,||\omega||_\infty|\gamma'(t)|
\end{equation}
is defined for almost every $t\in [0,L]$ (see \cite[Theorem 9A, p.\ 303]{wh:book} and subsection \ref{intrinsic}). In particular, we can define, for $a,b\in X$,
\begin{equation}\label{rhod}
d_\rho(a,b):=\inf \int_0^{L} (\sum_{i=1}^n\langle \rho_i(t),\gamma'(t)\rangle^2)^{1/2}\,dt\,,
\end{equation}
where the infimum is taken over all Lipschitz embeddings $\gamma:[0,L]\to X$ such that $\gamma(0)=a$ and
$\gamma(1)=b$.

To see that $d_\rho$ defines a metric, observe first that $d_\rho(a,b)<\infty$ because $X$ is rectifiably connected (see \ref{sobolev} and \cite[Section 17]{ch:lipschitz}). It remains to show that
$d_\rho(a,b)>0$ if $a\ne b$. Fix a Lipschitz embedding $\gamma:[0,L]\to X$. We
write $\gamma$ also for the image $\gamma([0,L])$. Now  (see \cite[3.2.19]{fed:gmt})
 $\cH_1$ almost every point $q\in \gamma$ is  tangent to $\gamma$ and satisfies
\begin{equation}\label{tanu}
\frac{\cH_1(\gamma\cap B(q, r))}{2r}\le 2\,
\end{equation}
for every $0<r<r(q)$, where $B(q, 2r(q))\subset W$.
Next, fix $x\in B(q, r(q))\cap\gamma$ such that the open subarc $\gamma(q,x)$ of $\gamma$ from $q$ to $x$ lies in $B(q, |q-x|)$.
Then $\gamma(q,x)$ is naturally oriented and can be viewed as a flat 1-chain
in $W$.
Consider the map $F_q:B(q, 2r(q))\to\rn$ given in \eqref{omap}.
By  \cite[4.2.9, 4.2.10, and 4.5.14]{fed:gmt}, we can find a flat 2-chain $S$ in $B(q, 2r(q))$
such that
$\partial S=\gamma(q,x)-[q,x]$ and that \begin{equation}\label{eka2}
4\pi|S|\le |\gamma(q,x)-[q,x]|^2\,
\end{equation}
(cf.\ the proof of Lemma \ref{flgconv}). Therefore, by  the Stokes
formula \eqref{w2}, and by  \eqref{eka2}, \eqref{tanu}, we obtain
\begin{equation}\label{eka3}
|\int_{\gamma(q,x)}\rho-F_q(x)|\,=\,|\int_{\gamma(q,x)}\rho-\int_{[q,x]}\rho\,|\,=\,|\int_{S}d\rho|\,\le\,\frac{25}{4\pi}||d\rho||_{\infty}|q-x|^2\,.
\end{equation}
On the other hand, it follows from the properties of BLD-mappings (see \cite[Proposition 4.13 and Theorem 6.8]{hr:duke} or \cite[(4.4)]{hs:duke})
that
\begin{equation}\label{eka4}
c\,|q-x|\,\le\,|F_q(q)-F_q(x)|=|F_q(x)|\,,
\end{equation}
where $c>0$ depends only on the local data near $\gamma$. (This requires that we choose $r(q)$ initially small enough so that we have uniform Ahlfors regularity and local linear contractibility for
points and radii in question, and moreover such that $F_q|B(q, 2r(q))\cap X$ is a BLD-mapping with uniform data; see \ref{cartBLD}.)

Now let $x\in B(q, r(q))\cap\gamma$ be arbitrary. Denote by $x'$ the point on $\gamma$ between $q$ and $x$ such that $|q-x|=|q-x'|$ and that the open subarc $\gamma(q,x')$ from $q$ to $x'$ lies in $B(q, |q-x'|)$. We deduce from the preceding, especially from \eqref{eka3} and \eqref{eka4},
that
$$
|q-x|=|q-x'|\,\le\,C\,|\int_{\gamma(q,x')}\rho\,|\,\le\,C\,\int_{\gamma(q,x)} (\sum_{i=1}^n\langle \rho_i,\gamma'\rangle^2)^{1/2}\,dt\,,
$$
where $C>0$ is independent of $q$ and $x$.
Since $\cH_1$ almost every point $q$ on $\gamma$ satisfies the preceding condition, we  easily obtain that
\begin{equation}\label{eka5}
|a-b|\,\le\,\cH_1(\gamma)
\,\le\,C\,\int_{\gamma} (\sum_{i=1}^n\langle \rho_i,\gamma'\rangle^2)^{1/2}\,dt\,.
\end{equation}
It follows that $d_\rho$ is a metric as desired.

Note that \eqref{eka5}, \eqref{eka0}, and the local quasiconvexity of $X$
give that the metric $d_\rho$ is locally bi-Lipschitzly equivalent to the Euclidean metric on $X$; that is,
for every point in $X$ there is a neighborhood and a constant $C\ge1$ such that
\begin{equation}\label{eka7}
C^{-1}|a-b|\,\le\,d_\rho(a,b)\,\le\,C\,|a-b|\,
\end{equation}
for all pair of points $a,b$ in the neighborhood.

\begin{proof}[Proof of Theorem \ref{sbility2}]
As mentioned earlier, the proof is similar to that in \cite{sul:bld}; only our framework is  slightly different.
Let $d_\rho$ be the metric on $X$ given  in \eqref{rhod}. By the Sobolev assumption, Theorem \ref{main2} gives that
the map $F_p:B(p, r)\cap X\to\rn$ given in \eqref{main2k} is bi-Lipschitz for every $p\in X$ and $r<r(p)$.
Here it is immaterial whether we use the Euclidean metric on $X$ or the metric $d_\rho$ (see the remark
just before the proof). On the other hand, the asymptotic estimate in \eqref{asymp}
implies that   $F_p$ is close to an isometry with respect to $d_\rho$ for points near $p$.

To establish the last claim, fix $x,y\in B(p,r)$ for $r>0$ small, and let $\gamma_{xy}$ be any rectifiable arc
joining the two points; by \eqref{eka5} and \eqref{eka7}, we may assume that $\gamma\subset B(p, r(p))$ and that $\cH_1(\gamma)\le C|x-y|$ for some constant $C\ge1$ independent of $x$ and $y$. As earlier, we  think of $\gamma$ as a flat 1-chain, so that
\begin{align*}
|F_p(x)-F_q(y)|&=|\int_{\gamma}dF_p|\,
\le\,C\,|x-y|\,||\rho-dF_p||_{\infty. B(p,r)}+|\int_\gamma\rho|\,\\
&\le\,C\,d_\rho(x,y)\,r\,||d\rho||_{\infty, B(p,r)}+
\int_{\gamma} (\sum_{i=1}^n\langle \rho_i,\gamma'\rangle^2)^{1/2}\,dt\,,
\end{align*}
where $C>0$ depends only on the local data.
Thus,
\begin{equation}\label{eka8}
|F_p(x)-F_p(y)|\,\le\,(1+\epsilon(r))\,d_\rho(x,y)\,,
\end{equation}
where $\epsilon(r)\to 0$ as $r\to 0$. Similarly we obtain
\begin{align*}
\int_{\gamma} (\sum_{i=1}^n\langle \rho_i,\gamma'\rangle^2)^{1/2}\,dt\,&\le\,
\int_{\gamma} (\sum_{i=1}^n(\langle \rho_i-d(F_p)_i,\gamma'\rangle^2)^{1/2}\,dt\,+\int_{\gamma}( \sum_{i=1}^n\langle d(F_p)_i,\gamma'\rangle^2)^{1/2}\,dt\,\\
&\le\,C\,|x-y|\,||\rho-dF_p||_{\infty, B(p,r)}+\cH_1(F_p\circ\gamma)\\
&\le\,C\,d_\rho(x,y)\,r\,||d\rho||_{\infty, B(p,r)}+\cH_1(F_p\circ\gamma)\,,
\end{align*}
whence
\begin{equation}\label{eka9}
(1-\epsilon(r))\,d_\rho(x,y)\,\le\,|F_p(x)-F_p(y)|\,,
\end{equation}
where $\epsilon(r)\to 0$ as $r\to 0$.

In conclusion, we have shown that $X$ can be covered by open sets such that both \eqref{eka8} and \eqref{eka9} hold for every pair of points in the neighborhood. By the results in \cite{shi1:smooth},
 \cite{gm:smooth}, we obtain that $X$ has a compatible smooth structure as desired.

The proof of Theorem \ref{sbility2} is complete.
\end{proof}

\subsection{Remarks}
(a) The Sobolev condition for $\rho$ in Theorem \ref{sbility2} can be replaced by the VMO-condition in \eqref{uusi}
as explained in \ref{vmormks}.
More generally, it suffices to assume that
 $\rho$ is such that the local BLD-mappings $F_q$ as in \eqref{omap} are always local homeomorphisms.

(b) It is not sufficient in Theorem \ref{sbility2} to only assume the existence of a globally defined Cartan-Whitney presentation.
Namely, for each $n\ge3$ there are homology $n$-manifolds $X$ that satisfy all the hypotheses in Theorem \ref{sbility2} except the Sobolev regularity for $\rho$, but that are not even topological manifolds, cf.\ \cite[9.1]{hr:duke}, \cite[2.4]{hs:duke}

(c) As a further example,  let $X\subset \R^N$ be a compact $n$-dimensional polyhedron that is also an orientable PL-manifold (with respect to the given  triangulation), but not smoothable; such polyhedra exist for many $n$  \cite{milnor:icm62}. Then we can map $X$ onto the standard $n$-sphere by a BLD-mapping and pull back  the standard coframe to obtain a global Cartan-Whitney presentation on $X$ \cite[2.4]{hs:duke}. (Note that the mapping has a Lipschitz extension to an open neighborhood of $X$ in $\R^N$.) We suspect that by using an example of this kind, one could demonstrate the sharpness of the Sobolev condition in Theorem \ref{sbility2}.

We do not know if analogous examples exist in dimensions where every PL-manifold is smoothable, e.g.\ in dimension $n=4$ \cite{whi:trans}.
In particular, we recall that  it is not known whether every Lipschitz 4-manifold is smoothable.

\

\subsection{Sullivan's cotangent structures}
Theorem \ref{sbility} is a essentially a direct consequence of Theorem \ref{sbility2}.
For completeness, we review the concept of a cotangent structure from \cite{sul:bld}
and include the proof.

Let $X$ be a Lipschitz $n$-manifold as defined in subsection  \ref{basic}, and let $\{(U_\alpha, \vp_\alpha)\}$ be a Lipschitz atlas. Thus,  $\{U_\alpha\}$ is an open cover
of $X$ and
 $$
 \vp_\alpha:U_\alpha\to \vp_\alpha(U_\alpha)=:\Va\subset \R^n
 $$
 are homeomorphisms such that the transition maps
 $$
\varphi_{\alpha\beta}=\varphi_\beta\circ\varphi_\alpha^{-1}:\varphi_\alpha(U_\alpha\cap U_\beta)\to \varphi_\beta(U_\alpha\cap U_\beta)\,
$$
are bi-Lipschitz.
A {\it Whitney} or {\it flat} $k$-{\it form} $\omega$ on $X$ is by definition a collection $\{\omega_\al\}$ of flat $k$-forms $\omega_\al$ defined in $V_\al$ such that
 the compatibility conditions
$$
\varphi_{\alpha\beta}^*\omega_\beta=\omega_\alpha
$$
hold in $\varphi_\alpha(U_\alpha\cap U_\beta)\subset V_\alpha$.
(Compare \cite{tel:ihes}.)
Flat $k$-forms on $X$  obviously form a vector space, that is moreover a module over the ring Lip$(X)$
of (locally) Lipschitz functions on $X$.

It makes sense to talk about Lipschitz vector bundles over $X$; the transition functions
are Lipschitz maps into a linear group $GL(\cdot,\R)$.
By using a Lipschitz partition of unity,
one can construct metrics, or Lipschitzly varying inner products, in Lipschitz vector bundles over $X$.

We say that $X$ admits a {\it cotangent structure} if there exists pair $(E,\iota)$, where $E$ is an oriented Lipschitz vector bundle of rank $n$ over $X$ and $\iota$ is a module map over Lip$(X)$ from Lipschitz sections of $E$ to flat 1-forms on $X$ satisfying the following:

 (*)
 If $s_1,\dots, s_n:X\to E$ are Lipschitz sections such that the induced section $s_1\wedge\dots\wedge s_n$ determines the chosen orientation on the bundle $\wedge_nE$, then for every $\al$  the flat $n$-form
 $
 \iota(s_1)_\alpha\wedge\dots\wedge \iota(s_n)_\alpha=\tau_\al\,dx_1\wedge\dots\wedge dx_n
 $
 in $V_\al$ satisfies
 \begin{equation}\label{alaraja}
 {\text{ess inf}}_{K} \  \tau_\alpha>0
\end{equation}
for every compact subset $K$ of $V_\alpha$. Here ess inf refers
to Lebesgue $n$-measure.

The above description of a cotangent structure is a bit simpler than that in \cite{sul:bld}; in particular, we make no mention about local (torsion free) connections.
This is possible because of the proof for the existence of local BLD-mappings in \cite[4.6]{hs:duke}. On the other hand, we consider oriented bundles for simplicity.

If $X$ is a smoothable, orientable $n$-manifold, then $X$ obviously admits a cotangent structure, namely we can put $E=T^*X$.

\begin{proof}[Proof of Theorem \ref{sbility}]
Let a Lipschitz  $n$-manifold $X$ admit a cotangent structure $(E,\iota)$. The Sobolev hypothesis in Theorem \ref{sbility} means that
  $\iota(s_i)_\alpha\in H^{1,2}_{loc}(V_\alpha)$ for every $\alpha$ and for every orthonormal frame $
  s_1,\dots, s_n$ of Lipschitz sections in $U_\alpha$. In particular, by Theorem \ref{main2}
  the local BLD-mappings considered by Sullivan in \cite[Theorem 1]{sul:bld} all have local degree equal to one. Theorem \ref{sbility} now follows as in the proof of Theorem 2 in \cite[p.\ 336]{sul:bld}.
(The proof there is analogous to the proof of Theorem \ref{sbility2}: one uses a global frame
to define a  metric on $X$ such that the local BLD-mappings are near isometries with respect to this metric.)
\end{proof}

\subsection{Remark} We stipulated in the preceding that the Sobolev hypothesis in Theorem \ref{sbility} means that $\iota(s_i)_\alpha\in H^{1,2}_{loc}(V_\alpha)$ for every $\alpha$ and for every orthonormal frame $s=(s_1,\dots, s_n)$ of Lipschitz sections in $U_\alpha$. This regularity hypothesis is independent of the frame, and is thus a property of the bundle only. Indeed, any two frames $s,s'$
differ  by a Lipschitz section $\sigma:X\to {\text{Aut}}E$, and $\iota$ is a module map over Lip$(X)$. Therefore, the mapping from $(\iota(s_1)_\alpha, \dots, \iota(s_n)_\alpha)$ to $(\iota(s_i')_\alpha, \dots, \iota(s_n')_\alpha)$ preserves the first order Sobolev spaces.

\bibliography{biblio}

\end{document}